\newif\ifpreprint
\newif\ifinforms

\preprinttrue %

\ifpreprint
\documentclass[12pt]{article}
\pdfoutput=1

\usepackage{fullpage}

\usepackage{cancel}
\usepackage{tcolorbox}
\usepackage{adjustbox}  %
\usepackage{amsmath,amssymb,amsthm}

\allowdisplaybreaks
\usepackage[colorlinks=true, citecolor=cyan, linkcolor=cyan, urlcolor=cyan, 
pdfborder={0 0 0}  %
]{hyperref}

\makeatletter
\renewcommand{\paragraph}[1]{\textbf{#1}\;}
\makeatother

\usepackage[round,authoryear]{natbib}
\fi

\newcommand{\myack}{I.\ Wang and B.\ Stellato are supported by the NSF CAREER Award ECCS-2239771. I.\ Wang is also supported by the Princeton Wallace Memorial Fellowship. B.\ Stellato is also supported by the ONR YIP Award N000142512147.}
\newcommand{\mythanks}{We thank Marta Fochesato for her contributions to this work, and Bart Van Parys for helpful discussions. The authors are pleased to acknowledge that the work reported on in this paper was substantially performed using Princeton University's Research Computing resources.}

\ifinforms
\documentclass[ijoc,nonblindrev]{informs4}
\OneAndAHalfSpacedXII

\usepackage{natbib}     %
\providecommand{\newblock}{} %
\usepackage{csvsimple}	%
\usepackage{booktabs}   %
\let\oldtabular\tabular
\let\endoldtabular\endtabular
\renewenvironment{tabular}{\oldtabular}{\endoldtabular}
\usepackage{adjustbox}  %

\TheoremsNumberedThrough
\EquationsNumberedThrough

\RUNTITLE{Fast Online DRO via Data Compression}
\RUNAUTHOR{Wang and Stellato}

\fi %

\ifinforms
\usepackage{xr}
\fi

\usepackage[font=small,
skip=5pt,
belowskip=0pt
]{caption} %

\usepackage{algorithm}%
\usepackage[noend]{algpseudocode} %

\usepackage{subfig}
\usepackage{verbatim}
\usepackage{graphicx}
\usepackage{mathtools}
\usepackage{enumitem}
\setitemize{leftmargin=3em}
\usepackage{multirow}

\usepackage{csvsimple}	%
\usepackage{booktabs}   %
\AtBeginEnvironment{tabular}{\small}
\usepackage{adjustbox}  %
\usepackage{appendix}
\usepackage{multicol}

\usepackage[acronyms]{glossaries}   %
\newacronym{SAA}{SAA}{sample average approximation}
\newacronym{LP}{LP}{linear program}
\newacronym{QP}{QP}{quadratic program}
\newacronym{MIQP}{MIQP}{mixed-integer quadratic program}
\newacronym{MIP}{MIP}{mixed-integer program}
\newacronym{MILP}{MILP}{mixed-integer linear program}
\newacronym{MINLP}{MINLP}{mixed-integer nonlinear program}
\newacronym{sBB}{sBB}{spatial branch and bound}
\newacronym{NLP}{NLP}{nonlinear program}
\newacronym{PWA}{PWA}{piecewise affine}
\newacronym{SVM}{SVM}{support vector machines}
\newacronym{ReLU}{ReLU}{rectified linear unit}
\newacronym{CPU}{CPU}{central processing unit}
\newacronym{GPU}{GPU}{graphics processing unit}
\newacronym{MPC}{MPC}{model predictive control}
\newacronym{ADMM}{ADMM}{alternating direction method of multipliers}
\newacronym{ADP}{ADP}{approximate dynamic programming}
\newacronym{FPGA}{FPGA}{field-programmable gate array}
\newacronym{DRO}{DRO}{distributionally robust optimization}
\newacronym{RO}{RO}{robust optimization}
\newacronym{MRO}{MRO}{mean robust optimization}
\newacronym{SO}{SO}{stochastic optimization}

\makeglossaries

\definecolor{purple}{RGB}{128,0,128}

\newcommand{\eg}{{\it e.g.}}
\newcommand{\ie}{{\it i.e.}}

\newcommand{\ones}{\mathbf 1}
\newcommand{\reals}{{\bf R}}

\newcommand{\prob}{{\mathbf P}}

\newcommand{\Expect}{{\mbox{\bf E}}}

\newcommand{\dom}{\mathop{\bf dom}_u} %

\newcommand{\define}{=}

\newcommand{\CVaR}{\mathop{\bf CVaR}}

\newcommand{\supp}{S}
\newcommand{\unc}{u}
\newcommand{\eps}{\varepsilon^K}
\newcommand{\epsnom}{\varepsilon}
\newcommand{\pwass}{p}  %

\newcommand{\gap}{\mathop{\bf gap}}
\ifpreprint
\newcommand{\argmin}{\mathop{\rm argmin}}
\newcommand{\argmax}{\mathop{\rm argmax}}
\newtheorem{assumption}{Assumption}[section]
\newtheorem{theorem}{Theorem}[section]  %
\newtheorem{lemma}{Lemma}[section]  %
\newtheorem{corollary}{Corollary}[section]  %
\newtheorem{proposition}{Proposition}[section]
\newtheorem{definition}{Definition}[section]
\newtheorem{example}{Example}[section]

\fi

\ifpreprint

\else

\fi 

\newcommand{\paperabstract}{%
We propose an online data compression approach for efficiently solving Wasserstein \gls{DRO} problems with streaming data.
Our method constructs adaptive ambiguity sets around a compressed distribution obtained by online clustering, so the problem size stays fixed as data accumulate.
We construct a dynamic regret bound with respect to a one-step-ahead non-compressed \gls{DRO} oracle, and establish online clustering conditions such that, with high probability, the regret converges sublinearly to a clustering-discrepancy-based performance gap. This gap is defined in terms of the discrepancies between the true and compressed distributions, so that by varying the number of clusters, our method trades off robustness against computational effort.
We additionally provide fast subgradient-based updates that replace direct solutions of both the full and compressed problems, and extend the regret analysis to this setting.
Numerical experiments on portfolio optimization and sparse support vector machine problems, including mixed-integer formulations, show over an order of magnitude reduction in cumulative computation time, even compared to non-robust methods, with minimal loss in solution quality.}

\newcommand{\paperkeywords}{Distributionally robust optimization, Data-driven optimization, Streaming data, Online clustering, Dynamic regret}

\ifpreprint
\title{Fast Online Distributionally Robust Optimization \\via Data Compression }
\author{\large Irina Wang and Bartolomeo Stellato}
\fi

\ifinforms
\TITLE{Fast Online Distributionally Robust Optimization via Data Compression}
\ARTICLEAUTHORS{%
\AUTHOR{Irina Wang}
\AFF{Department of Operations Research and Financial Engineering, Princeton University, New Jersey, USA, \EMAIL{iywang@princeton.edu}}
\AUTHOR{Bartolomeo Stellato}
\AFF{Department of Operations Research and Financial Engineering, Princeton University, New Jersey, USA, \EMAIL{bstellato@princeton.edu}}
}
\ABSTRACT{\paperabstract}
\FUNDING{\myack}
\fi

\begin{document}

\maketitle

\ifpreprint
\begin{abstract}
\paperabstract
\end{abstract}
\fi

\glsresetall

\section{Introduction}
Many decision-making problems in engineering, operations research, and computer science involve solving optimization problems affected by uncertain parameters, and
\gls{DRO} has emerged as a disciplined framework to robustify the decision-making process against distributional misspecification~\citep{dro_survey,kuhn_dro_book}.
Instead of assuming a known distribution, \gls{DRO} considers an \emph{ambiguity set} of possible distributions and optimizes for the worst-case expected cost under all distributions within this set.
This approach offers out-of-sample performance guarantees with a finite number of data samples, overcoming the {\it optimizer's curse} of~\gls{SAA}.

Since its appearance, significant effort has been devoted to constructing ambiguity sets primarily in static settings, using either structural assumptions on the distribution (\eg, moments~\citep{moment,moment_1}) or training data directly~\citep{kuhn2019wasserstein,gao2016distributionally}.
However, many real-time settings, such as healthcare resource allocation, air traffic control, and energy trading, are dynamic and new data becomes available sequentially over time.
Despite its potential,~\gls{DRO} remains underexplored in such scenarios.
Most~\gls{DRO} approaches construct a static ambiguity set, and model time progression through {\it adjustable wait-and-see variables}~\citep{adro2,adro,robustadaptopt}.
There are limited works focusing on how to adapt the ambiguity sets dynamically over time, and each targets a specific problem class: distributions supported on a finite set of scenarios~\citep{drotime}, continuous variables with locally strongly concave objectives~\citep{assim}, or continuous variables with first-order updates tailored to Wasserstein-1 sets and losses with strictly sublinear growth in the uncertainty~\citep{online_dro_s}.
In particular, streaming data settings can be challenging for the popular Wasserstein-based \gls{DRO}~\citep{mohajerin2018data,kuhn2019wasserstein,gao2016distributionally}, where the ambiguity set is defined as a Wasserstein ball around the empirical distribution of a dataset.
As new datapoints arrive, the problem dimension grows, increasing computational complexity~\citep{kuhn2019wasserstein}.
While a larger dataset improves confidence in the empirical measure as an approximation of the true distribution, the reduction in the ambiguity set’s radius is not enough to offset this complexity; for mixed-integer problems, even \gls{SAA} is difficult to solve with large amounts of data, without any added complexity from ambiguity. This leads to computational bottlenecks that are especially limiting for real-time applications.

\subsection{Our contributions}%
\label{sub:our_contributions}
We propose an {\rm online data compression approach} for efficiently solving Wasserstein~\gls{DRO} problems with streaming data.
Our key contributions are:
\begin{itemize}
    \item {\it Online \gls{DRO} framework with data compression:}
    We develop a generic framework for solving Wasserstein~\gls{DRO} problems with streaming data, where {\it adaptive} ambiguity sets are formulated as gradually shrinking Wasserstein balls, of any order, around a clustered empirical distribution.
    Since the compressed reference distribution is supported on at most a fixed number of clusters, the problem size does not grow with the data stream. %
    At every time step, the decision-maker can either solve the compressed problem directly, or apply fast approximate subgradient updates.
    We show that our method is compatible with general clustering procedures, and provide examples fast and memory-efficient online clustering algorithms.
    \item {\it Online regret analysis:}
    We construct a dynamic regret bound comparing our online solutions to a one-step-ahead nominal \gls{DRO} oracle, \ie, with access to the subsequent realization of the uncertainty.
    We show that, with high probability, the regret converges sublinearly to a function of the discrepancy between the true and compressed distributions. 
    This quantifies the impact of data compression, and highlights the {\it trade-off between computational effort and optimality.}
    \item {\it Fast subgradient subroutines:}
    For maximum-of-affine losses, we decompose the worst-case expectation into per-cluster oracles with closed-form or efficiently computable solutions, under various Wasserstein orders, ground norms, and supports, together with per-iteration complexity and tolerance guarantees.
    These subroutines replace direct solves of the compressed and full \gls{DRO} problems, and we extend the dynamic regret analysis to this setting.
    \item {\it Computational gains in numerical experiments:}
    We demonstrate the efficiency of our approach on portfolio optimization and sparse (cardinality-constrained) support vector machine problems with streaming data.
    Our results show significant computational savings, even compared to~\gls{SAA}, with minimal loss in solution quality.
\end{itemize}

\subsection{Related work}%
\label{sub:related_work}

\paragraph{Distributionally robust optimization}
\Gls{DRO} has been extensively explored in recent years, with successful applications in machine learning \citep{shafieezadeh2015distributionally}, among many other areas.
We focus on discrepancy-based ambiguity sets, defined as a ball around a data-driven nominal distribution, where the distance, expressed \eg\ via the $\phi$-divergence \citep{bayraksan2015data} or the Wasserstein distance \citep{mohajerin2018data,kannan2024residuals}, signifies the trust in the data at hand.
Due to its expressivity and favorable statistical properties, we focus on Wasserstein \gls{DRO}.

\paragraph{Online learning and \gls{DRO}}
Online learning provides algorithms for solving repetitive problems over time, and has recently been applied to robust optimization \citep{ho2018online,chen2017robust} and to \gls{DRO}: \cite{namkoong2016stochastic} and \cite{qi2021online} propose online stochastic methods for $\phi$-divergence ambiguity sets or KL-divergence regularization.
We focus here on a different problem, where the same \gls{DRO} problem is solved repeatedly over time with growing knowledge of the uncertainty distribution; our work is thus closest to the few approaches that adapt the ambiguity set to the samples collected progressively.
\cite{drotime} update the ambiguity set through projected gradient descent, but consider distributions supported on a fixed and finite set of scenarios, a setting in which the problem size does not grow with the data stream.
\cite{assim} incorporate streaming data through assimilation, but require continuous variables and strongly concave objectives.
Most recently, \cite{online_dro_s} solve the online~\gls{DRO} problem through first-order updates, with efficient subroutines for losses with strictly sublinear growth in the uncertainty, for Wasserstein order one; this is similar in spirit to our subgradient approach, which instead covers max-of-affine losses for any order, so the two works are complementary.
Unlike these purely first-order approaches, our framework also supports direct solves, and in particular mixed-integer formulations, which are our main computational focus.

\paragraph{Scenario reduction and \gls{DRO}}
Since the size and computational complexity of data-driven optimization problems generally increase with the number of samples, scenario reduction techniques have been developed to reduce the number of scenarios while retaining an accurate representation of the underlying uncertainty.
Common approaches minimize the distance between the reduced and the original distribution \citep{dupavcova2003scenario}, or aggregate scenarios via clustering \citep{hewitt2022decision}, moment matching \citep{mehrotra2013generating}, objective approximation \citep{zenios1993constructing}, nested distances \citep{horejvsova2020evaluation}, and, more recently, decision-focused metrics that involve the loss function itself \citep{bertsimas2023optimization, zhang2023optimized}.
Our previous work~\citep{wang2022mean} uses clustering in static, finite-sample \gls{DRO} problems as a tool to bridge robust and distributionally robust optimization.
Closely related, \cite{aigner2025scenario} embed scenario reduction into a \gls{DRO} framework with suboptimality bounds, albeit limited to monotonically homogeneous uncertain objectives with strictly positive uncertainty; both works, however, assume data are available a priori, while we focus on {\it online} scenario reduction for \gls{DRO} problems with streaming data.

\paragraph{Clustering with streaming data}
Data stream clustering methods can be broadly classified into partition-based algorithms, which group data using distance-based similarity metrics, density-based algorithms such as DenStream \citep{cao2006density}, and hierarchical algorithms \citep{murtagh2012algorithms}; see the survey~\citep{streaming} for details.
Among the first class, the most popular ones are incremental $k$-means \citep{aaron2014dynamic} and CluStream \citep{clustream}, which we adapt to fit our framework, although our approach can incorporate any online clustering algorithm, with various memory and run-time trade-offs.

\subsection{Notation}
Let $\|\cdot\|$ denote an arbitrary norm on $\reals^d$, and let $\|\cdot\|_*$ denote its dual norm.  Let $\Xi_{\pwass}(S)$ be the space of all probability distributions $\mathbf{Q}$ defined on a support~$\supp \subseteq \reals^d$, with bounded $\pwass$-th moments, \ie,
$\Expect_\mathbf{Q}[\|u\|^{\pwass}] = \int_S\|u\|^{\pwass}\mathbf{Q}(du) < \infty$, where either $\pwass\geq 1$ or $\pwass = \infty$. 
Let $\mathbf{E}_{\mathbf{P}}[f(\unc,x)]$ denote the expectation of $f(\unc,x)$ with respect to $\unc \in\mathbf{P}$, and let $\dom f$ denote the domain of $f$ in the variable $u$. 
Let $[f]^*(z, x) \define \sup_{u\in \dom f} z^Tu - [f(u, x)] $ be the conjugate of $f$, and let $\sigma_\supp(z) \define\sup_{u \in  \supp}  z^Tu$ be the support function of $\supp$.
Let $q$ be the conjugate number of $\pwass$, satisfying $1/\pwass + 1/q = 1$, \ie, $q = \pwass/(\pwass - 1)$. Let $(a)_+ = \max(a,0)$. We write $\delta_u$ as the Dirac distribution concentrating unit mass on $u$, and $\Pi_\mathcal{X}$ as the projection operator onto a closed and convex set $\mathcal{X}$. 

\section{The online stochastic problem}
\label{sec:problem_intro}
We consider a \gls{SO} problem of the form
\begin{equation}
\label{eq:opt}
 H_\star  = \min_{x\in \mathcal{X}} \mathbf{E}_{\mathbf{P}}[f(\unc,x)],
\end{equation}
where $x \in \mathcal{X} \subseteq \mathbf{R}^n$ is the decision variable, $\mathcal{X}$ a compact set, $\unc \in \supp \subseteq \mathbf{R}^d$ the vector of uncertain parameters that is governed by some probability distribution $\mathbf{P}$, and $H_\star$ the optimal objective value. 
We assume that the support $\supp$ of $\prob$ lives within the domain of $f$ for the variable $u$, \ie, $\supp \subseteq \dom{f}$, and is closed and convex. 
The function $f:\mathcal{X}\times \supp \rightarrow (-\infty, \infty]$ is assumed to be of the form $f(u, x) = \max_{j\leq J}f_j(u,x)$, with each $f_j$ being proper, concave, and upper-semicontinuous in~$u$ for all $x$, and convex in $x$.
We make the following two regularity assumptions on $f_j$, where the existence of global constants independent of $x$ is implied by the compactness of $\mathcal{X}$.

\begin{assumption}
\label{ass:general_ass_lip} For all $x \in \mathcal{X}$, the functions $f_j$ are Lipschitz continuous in $u$ with constants~$L_j$, \ie, $\left |f_j(v,x) - f_j(u,x)\right| \leq L_j\| u - v\|, \; \forall u,v \in \dom f$. Let $L = \max_{j \leq J} L_j$.
\end{assumption}

\begin{assumption}
\label{ass:general_ass_smooth} For all $x \in \mathcal{X}$, the functions $f_j$ are smooth (w.r.t.\ the $\ell_2$ norm) in $u$ with constants~$M_j$, \ie, $\left\|\nabla f_j(v,x) - \nabla f_j(u,x)\right\|_2 \leq M_j\| u - v\|_2,\; \forall u,v \in \dom f$. Let $M = \max_{j \leq J} M_j$.
\end{assumption}
Analogously to~\citep{wang2022mean}, we assume the following on the domain of $f$. 
\begin{assumption}
	\label{assp:dom}
	The domain $\dom{f}$ is $\reals^d$. Otherwise, $f$ is either element-wise monotonically increasing in $u$ and only has a (potentially) lower-bounded domain, or element-wise monotonically decreasing in $u$ and only has a (potentially) upper-bounded domain. In all cases, the domain $\dom{f}$ is independent of $x$.
\end{assumption}
Assumptions~\ref{ass:general_ass_lip} and~\ref{assp:dom} hold throughout the paper, and we do not restate them in subsequent results.
The smoothness condition of Assumption~\ref{ass:general_ass_smooth} is instead invoked only where explicitly stated.

\paragraph{Streaming data}
We assume the distribution $\prob$ to be unknown, so the value $H_\star$ cannot be computed directly.
Instead, we have access to a {\it streaming dataset} of i.i.d.\ realizations of $u$: starting from an initial dataset $\mathcal{D}_0 = \{\hat{u}^i\}_{i=1}^{n_0}$ with $n_0 \geq 1$ datapoints, one realization is disclosed (w.l.o.g.) at the end of each time step $t \geq 1$, so that $\mathcal{D}_t = \{\hat{\unc}^{i}\}_{i=1}^{n_t}$ with $n_t = n_0 + t$.
We define the corresponding empirical distribution $\hat{\prob}_t = (1/{n_t})\sum_{i=1}^{n_t} \delta_{\hat{\unc}^{i}}$.

\subsection{\Gls{DRO} over time}
\label{sec:over_time}
\label{subsec:nom_perf}

 The~\gls{DRO} approach under the Wasserstein metric constructs an ambiguity set of order $\pwass$ around $\hat{\prob}_t$ of radius $\epsnom_t \geq 0$, \ie,
$\mathcal{P}_t = \left\{\mathbf{Q} \in \Xi_{\pwass}({\supp}) \mid W_\pwass(\hat{\prob}_t, \mathbf{Q}) \leq \epsnom_t\right\},$
where the Wasserstein type-$\pwass$ distance between two probability distributions $\mathbf{Q}_1, \mathbf{Q}_2 \in \Xi_{\pwass}(\supp)$ is defined by
$W_{\pwass}(\mathbf{Q}_1, \mathbf{Q}_2)=\inf _{\pi \in \Pi(\mathbf{Q}_1, \mathbf{Q}_2)} (\int_{\supp \times \supp} \|\unc_1-\unc_2\|^{\pwass} \mathrm{d} \pi(\unc_1, \unc_2))^{1/{\pwass}}$, with~$\Pi(\mathbf{Q}_1, \mathbf{Q}_2)$ the set of all couplings, \ie, probability distributions over $\supp \times \supp$ with marginals $\mathbf{Q}_1$ and $\mathbf{Q}_2$.
We suppress the dependency of $\mathcal{P}_t$ on $\pwass$ for conciseness, and note that unless otherwise specified, $\pwass$ is of an arbitrary order and $\|\cdot\|$ is an arbitrary norm.
Intuitively, the ambiguity set~$\mathcal{P}_t$ contains all probability distributions in $\Xi_{\pwass}(\supp)$ that can be obtained by transporting probability mass from $\hat{\prob}_t$ with a transportation budget of at most $\epsnom_t$.
With this uncertainty modeling framework in place, we compute the~\gls{DRO} solution $x_t \in \mathcal{X}$ at time step~$t$ by solving the optimization problem
\begin{equation}\label{eq:DRO}
H_t = \min_{x \in \mathcal{X}} \max_{\mathbf{Q} \in \mathcal{P}_{t}}\mathbf{E}_{\mathbf{Q}}[f(\unc,x)],
\end{equation}
where we minimize the worst-case expected value of $f$.

As~$t$ increases, our confidence in $\hat{\prob}_t$ grows, allowing us to safely reduce the {\it nominal radius}~$\epsnom_t$; suitable choices of $\epsnom_t$ yield {\it finite sample performance guarantees} of the form
\begin{equation}
	\label{eq:prob_guarantees1}
	\prob^{t}\left( H_\star \leq \Expect_{\prob}[f(u, x_t)] \le H_t\right) \ge 1 - \beta_{t},
\end{equation}
where $\beta_{t} > 0$ is a specified time-varying probability of constraint violation, and $\prob^{t}$ is the product distribution of the dataset $\mathcal{D}_{t}$.
We choose the radius according to the classic literature~\citep{fournier2015rate,kuhn2019wasserstein}, which, for finite $p$, requires the following light-tailed assumption.  Note that it is satisfied by sub-Gaussian and bounded distributions, and holds throughout the paper.
\begin{assumption}[Light-tailed distribution]
\label{ass:light_tail}
    There exists an exponent $\alpha > 1$, such that
    $R = \Expect_\prob[\exp(\|u\|^{\alpha})] = \int_S \exp(\|u\|^{\alpha})\prob(du) < \infty.$
\end{assumption}
Then, measure concentration results for Wasserstein-$\pwass$ distances imply the following.
\begin{theorem}[Guarantees using Wasserstein measure concentration]
    \label{thm:calc_eps}
    Under Assumption~\ref{ass:light_tail} for the case $\pwass \geq 1$, and under analogous assumptions for $\pwass =\infty$~\citep{wass-p-guarantee}, for $n_t \rightarrow \infty$ we can choose confidence levels $\beta_t \in (0,1)$ and corresponding radii $\epsnom_t=O((\log(\beta^{-1}_t)/{n_t})^{\min(\pwass/d,1/2)})$ such that $\sum_{t=1}^\infty \beta_t < \infty$,
    $\prob^{t}(W_{\pwass}(\prob,\hat{\prob}_t) \geq \epsnom_t) \leq \beta_t,$ and $\lim_{t \rightarrow \infty} \epsnom_t = 0.$ The~\gls{DRO} problem~\eqref{eq:DRO} solved with these radii gives the finite-sample guarantee~\eqref{eq:prob_guarantees1}.
\end{theorem}

Larger datasets thus reduce the distributional ambiguity: Theorem~\ref{thm:calc_eps} gives $\mathbf{P}^\infty \{ \lim_{t \rightarrow \infty} W_\pwass(\mathbf{P}, \hat{\mathbf{P}}_t) =0 \} = 1$, which in turn implies $\mathbf{P}^\infty \{\lim_{t \rightarrow \infty}H_t \downarrow H_\star\}=1$~\citep[Theorem 3.7]{mohajerin2018data}.
This statistical benefit, however, comes at a steep computational price: the radius decays only as $O(n_t^{-\pwass/d})$, suffering from the curse of dimensionality, while the reformulation of the \gls{DRO} problem~\eqref{eq:DRO} adds a constraint for each datapoint in $\hat{\mathbf{P}}_t$, so the problem size grows without bound as $t$ increases.
We are thus motivated to introduce a new tradeoff, where some performance is sacrificed for faster computations.

\subsection{Our online algorithm}
\label{subsec:ours}
We propose a data-compression framework for fast online \gls{DRO}.
Our goal is to solve an adjusted \gls{DRO} problem over time, where $K$ is the chosen number of data clusters, and
\begin{equation}\label{eq:dro:time}
    \begin{aligned}
x^K_t &= \argmin_{x \in \mathcal{X}} \max_{\mathbf{Q} \in \mathcal{P}^K_{t}}\mathbf{E}_{\unc \sim \mathbf{Q}}[f(\unc,x)],\quad t =0,1,\dots,\\
    \end{aligned}
\end{equation}
where $\mathcal{P}^K_{t} = \{\mathbf{Q} \in \Xi_{\pwass}({\supp}) \mid W_\pwass(\hat{\prob}^K_t, \mathbf{Q}) \leq \eps_t\}$ is an {\it adaptive ambiguity set} with radius $\eps_{t} \geq 0$, constructed from the data $\mathcal{D}_{t}$ observed so far. From here onward, we refer to the solution $x^K_t$ as the {\it online} solution, and the solution $x_t$ from~\eqref{eq:DRO} as the {\it nominal}~\gls{DRO} solution.

Unlike the standard \gls{DRO} ambiguity set, the adaptive ambiguity set is {\it compressed}: its reference distribution is supported on at most $K \leq n_t$ points, obtained by {\it clustering} the datapoints into sets $C^k_t$, $k = 1, \dots,K_t,$ where $K_t \leq K$.
Each cluster $C^k_t$, with $n^k_t = |C^k_t|$ datapoints, has weight $\theta_t^k = n^k_t/n_t$ and mean $\bar{u}^k_t = (1/n^k_t) \sum_{\hat{u} \in C^k_t} \hat{u}$, which characterize the clustered empirical distribution $\hat{\mathbf{P}}^K_t = \sum_{k=1}^{K_t}\theta_t^k\delta_{\bar{u}_t^k}$.
In Section~\ref{sec:shared_cluster}, we provide additional details on our clustering requirements, and in Section~\ref{sec:algos}, provide efficient online clustering methods.

We solve the online compressed~\gls{DRO} problem~\eqref{eq:dro:time} through the direct reformulation~\citep[Section 2.4]{wang2022mean} given below, to a tolerance $\delta$, which could be specified by modern solvers.
\begin{equation}
	 \label{eq:robustopt_p_max}
	\begin{array}{lll}
 \!\!H^K_t \!= \!\!\!&\underset{x\in\mathcal{X},\lambda_t\geq0,s_t,z_t,y_t}{\mbox{minimize}} & \sum_{k=1}^{K_t} \theta^k_{t} s^k_t\\
&\mbox{subject to}  & [-f_j]^*(z^{jk}_t-y^{jk}_t , x) + \sigma_{{\supp}}(y^{jk}_t) - (z{^{jk}_t})^T\bar{\unc}^k_{t} + \dfrac{1}{q}\dfrac{\|z^{jk}_t\|_*^q}{(\lambda_t p)^{q-1}}   \\
	&	&\hspace{8em}   +\lambda_t \epsnom_{t}^p \leq s^k_t,\quad k = 1,\dots, K_t, \, j =1, \dots, J.\\	\end{array}\!\!\!\!\!\!\!\!\!\!\!\!\!
\end{equation}
In addition to the primal variables, we introduce dual and auxiliary variables $\lambda_t \in \reals$, $s^k_t \in \reals$, $z^{jk}_t \in \reals^{d}$, and $y^{jk}_t \in \reals^d$, for $k=1,\dots,K_t$ and $j=1,\dots,J$~\citep[Theorem 8]{kuhn2019wasserstein}.
This reformulation is for $\pwass \in [2,\infty)$; the cases for $\pwass = 1$ and $\pwass = \infty$ are given in~\citep{wang2022mean}. 

The computational complexity of~\eqref{eq:robustopt_p_max} scales with $K_t$ rather than with the total number of observed datapoints, significantly decreasing solution time compared to nominal~\gls{DRO}.
However, since even this problem can be expensive to solve at every timestep with low $\delta$, we can instead solve after every {\it batch} of data, or update the solution through a subgradient-based approach.

\subsubsection{Subgradient-based method}
\label{sec:subgrad}
For convex $\mathcal{X}$, we compute the worst-case expectation with respect to the current solution $x^K_{t-1}$, then take a projected subgradient step in $\mathcal{X}$ with some step size $\eta_t > 0$; see Algorithm~\ref{alg:hl}.
For non-convex $\mathcal{X}$, the same step applies only as a local-improvement heuristic; see Section~\ref{sub:sparse_svm} for an example.
\begin{algorithm}[h]
	\caption{Online Fast Compressed DRO Algorithm}
	\label{alg:hl}
	\begin{algorithmic}[1]
	 \State {\bf given} initial data $\mathcal{D}_0$, tol $\delta, \delta_{\text{approx}}$, cluster number $K$, sets $T_{\text{full-solve}}$ and $T_{\text{approx-solve}}$, step sizes $\eta_t$
\For{$t=0,1,\dots$}
     \State Use a subroutine to construct the clustered empirical distribution $\hat{\prob}^K_t$ \Comment{See Section~\ref{sec:algos}}
     \State Choose a radius for the ambiguity set $\mathcal{P}^K_{t}$ \Comment{See Sections~\ref{subsec:nom_perf} and~\ref{sec:exampless}}
       \If{$t  =  0~\text{or}~t\in T_{\text{full-solve}}$}
     \State $H^K_t, x^K_t \gets$ Solve~\eqref{eq:robustopt_p_max}, tol $\delta$ \Comment{See~\cite[Section 2.4]{wang2022mean}}
     \ElsIf{$t \in T_{\text{approx-solve}}$ and $\mathcal{X}$ is convex}\Comment{See Section~\ref{sec:subgrad_sec}}
     \State $ \hat{H}^K_t \gets \max_{\mathbf{Q} \in \mathcal{P}^K_{t}}\mathbf{E}_{\mathbf{Q}}[f(\unc,x^K_{t-1})]$, tol $\delta_{\text{approx}}$  \Comment{Solve for worst-case expectation}
     \State $g_t \in \partial_{x} \hat{H}^K_t  $ \Comment{Obtain a subgradient}
    \State $x^K_{t} \gets \Pi_{\mathcal{X}}(x^K_{t-1} - \eta_tg_t)$ \Comment{Update decision via projected subgradient method}
     \EndIf
     \State Observe a new realization of the uncertainty $\unc$, denoted as $\hat{\unc}$
  \EndFor
	\end{algorithmic}
\end{algorithm} 
The bottleneck of the subgradient-based method lies in step 8, finding the worst-case expectation at the fixed point $x^K_{t-1}$, which uses essentially the same reformulation as the full~\gls{DRO} problem~\citep{mohajerin2018data}, just without optimizing over $x$: it still introduces $O(KJ)$ auxiliary variables, costing $O((KJd)^{1.5})$ per solve (or $O((KJ(m+d))^{1.5})$ if the support $\supp$ is polyhedral with $m$ facets) via a generic interior-point method~\citep{nesterov1994interior}.
In Section~\ref{sec:subgrad_sec}, we instead construct specialized subroutines that solve this problem to tolerance $\delta_{\text{approx}}$ for maximum-of-affine $f$, for general $p \geq 1$ and $p=\infty$.
If the subroutines are faster than the clustering algorithm, they can even be run with $K = n_t$, bypassing the data-compression step.

\subsubsection{Dynamic regret}
In Section~\ref{sec:regret}, we analyze the {\it dynamic regret}
\begin{equation}
\label{regret}
R(T,K) = \frac{1}{T}\sum_{t=0}^{T-1} \left(\max_{\mathbf{Q} \in \mathcal{P}_{t+1}} \mathbf{E}_{\mathbf{Q}}[f(\unc,x^K_{t})] - \min_{x \in \mathcal{X}} \max_{\mathbf{Q} \in \mathcal{P}_{t+1}} \mathbf{E}_{\mathbf{Q}}[f(\unc,x)]\right),
\end{equation}
which compares the average performance of our online solution $x^K_{t}$ (solved to tolerance $\delta$) against that of a one-step-ahead oracle: the nominal \gls{DRO} solution computed with the updated, non-compressed ambiguity set $\mathcal{P}_{t+1}$.
This notion of regret is in line with the online learning literature~\citep{drotime,regret, regret2}.
We give deterministic finite-sample and asymptotic bounds, and additionally show sublinear convergence to the asymptotic bound. 
In Section~\ref{sec:subgrad_sec}, we extend the analysis to subgradient solves.

\section{Online clustering requirements}
\label{sec:shared_cluster}
\label{sec:general_intro}
The key to our approach is the adaptive ambiguity set $\mathcal{P}^K_t$, which is constructed using a {\it compressed} representation of the available data, $\hat{\prob}^K_t$.
Below, we formalize the requirements on generating $\hat{\prob}^K_t$, and note that the decision-maker has the flexibility to use {\it any} method to generate $\hat{\prob}^K_t$, provided it satisfies the given assumptions.
These requirements ensure the validity of the asymptotic regret bounds established later in this paper, and can be satisfied by any partitioning of the support set $\supp$.

The number of clusters $K$ reflects the available computational budget, as it caps the number of atoms in $\hat{\prob}^K_t$, and hence the size of the compressed problem~\eqref{eq:robustopt_p_max}.
For all methods and time periods $t$, we maintain a set of $K_t \leq K$ clusters $\{C_t^k\}_{k=1}^{K_t}$ and a corresponding set of supports $\mathcal{S}_t = \{S_t^k\}_{k=1}^{K_t}$, defined as follows.
\begin{definition}[Cluster support]
\label{def:cluster_supp}
    For each cluster $k$ and time $t$, the {\it cluster support} $S_t^k \subset S$ is defined as a region in $S$ such that any datapoint $\hat{u} \in \mathcal{D}_t$ that falls within this region is contained in cluster $C_t^k$.
    We say that the clusters $\{C_t^k\}_{k=1}^{K_t}$ are {\it induced} by the supports $\mathcal{S}_t$.
\end{definition}

\paragraph{Optimal clustering}
If the distribution $\prob$ is known, then for any~$K$, the optimal set of supports can be determined using $k$-centers clustering~\citep[Chapter 1]{quantization}.
\begin{definition}[$k$-centers]
    \label{def:k_center}
    The $k$-centers problems, with respect to a distribution $\mathbf{Q}$ with support $Z$, and orders $\pwass \geq 1$ (left), $\pwass=\infty$ with compact $Z$ (right), are
\begin{equation}
    \label{def:k_centers}
    \inf_{A \subset \reals^d, |A|\leq K}\Expect_\mathbf{Q}\left[\min_{a \in A}\|u - a\|^{\pwass}\right], ~\quad\text{and}\quad \inf_{A \subset \reals^d, |A|\leq K}\sup_{u \in Z}\min_{a \in A}\|u - a\|,
\end{equation}
where $A^\star$ is a set that achieves the infimum.
\end{definition}
Given optimal centers $A^\star$, each cluster support is a cell of the Voronoi diagram of $A^\star$,
$ S^k =\left\{u \in S\,\middle|\,\|u - a^k\| = \min_{a \in A^\star}\|u - a\| \right\},$
where each point $u$ is assigned to the region of its closest center $a^k$.
This clustering attains the {\it minimum quantization error}
$$w_K(\mathbf{Q}) = \left(\Expect_{\mathbf{Q}}\left[\min_{a \in A^\star}\|u - a\|^{\pwass}\right]\right)^{1/{\pwass}} ~~(\pwass \geq 1), \qquad w_K(\mathbf{Q}) = \sup_{u\in Z}\min_{a \in A^\star}\|u - a\| ~~(\pwass = \infty),$$
which is the smallest Wasserstein-$\pwass$ distance between $\mathbf{Q}$ and any distribution supported on at most $K$ atoms~\citep[Chapter 1]{quantization}.

Without knowledge of $\prob$, the optimal clustering with respect to the empirical distribution $\hat{\mathbf{P}}_t$ can be found in the same way.
However, the $k$-centers problem is NP-hard~\citep{k-centers}, and its size grows with $n_t$, conflicting with our goal of decreasing computational effort.
The goal of an online clustering algorithm is to cheaply approximate the $k$-centers solution.
Below, we give the assumptions the cluster supports must satisfy; Section~\ref{sec:algos} gives two algorithms.

\begin{assumption}
\label{ass:support}
The following holds for any time step $t$: the elements in the set $\mathcal{S}_t = \{S_t^k\}_{k=1}^{K_t}$ have pairwise disjoint interiors, \ie, $\mathrm{int}(S^k_t) \cap \mathrm{int}(S^{k'}_t) = \emptyset$ for $k \neq k'$; the cluster supports cover $S$, \ie, $S  = \cup_{k\leq K_t}S_t^k$; the cluster supports have nonzero measure, \ie, $\prob(S^k_t) > 0$ for $k=1,\dots,K_t$; and the cluster boundaries have measure 0.
Furthermore, there is some finite time $\tau < \infty$ such that $\mathcal{S}_t = \mathcal{S}_\tau$ for all $t \geq \tau$.
\end{assumption}
Assumption~\ref{ass:support} holds throughout the remainder of the paper.
By this assumption, $\hat{\mathbf{P}}^K_t$ represents an approximation of the  empirical distribution $\hat{\mathbf{P}}_t$, and converges to some distribution $\prob^K_\star$ as $t \to \infty$. 
\begin{lemma}[Clustering consistency]
\label{lem:cluster_converge}
The sequence of clustered empirical distributions~$\hat{\prob}^K_t$ converges~$\prob^\infty$-almost surely to a distribution $\prob^K_\star$. Furthermore, the distributions achieve almost sure convergence with respect to the Wasserstein metric, that is, for any $\pwass$ such that $\prob$ has finite $\pwass$-th moments,
    $\prob^{\infty} \left \{ \lim_{t \to \infty} W_{\pwass}(\prob^K_\star, \hat{\prob}^K_t) = 0\right\} = 1.$
\end{lemma}
\paragraph{Remark on $\tau$} We enforce the support-fixing time $\tau$ to ensure the convergence of the clustered distribution, which is required by our subsequent asymptotic results. 
Without a fixing time, convergence can only be ensured for specific algorithms: \cite{pmlr-v151-so22a} adjust the $k$-means center update rate to obtain asymptotic convergence, while \cite{pmlr-v54-tang17b} impose a uniform and deterministic update rate to further obtain a $O(1/t)$ convergence rate.
Both approaches fit in our framework, with the cluster supports set as the Voronoi cells of the obtained centers.
However, they require specific and difficult-to-compute constants, which motivates our more general support-fixing condition.

\section{Online performance analysis}
\label{sec:regret}
We are now ready for our regret analysis. 
In online settings, we want to robustify the solution against the {\it subsequent realization} of the uncertainty; the dynamic regret~\eqref{regret} therefore compares our online solution, with information up to time $t$, against the nominal \gls{DRO} oracle with information up to time $t+1$.
In order to obtain a bound on~\eqref{regret}, we first introduce a number of prerequisite lemmas. 
We begin by formalizing the performance difference between two general Wasserstein \gls{DRO} problems with potentially different centers and/or radii. 

\begin{lemma}[Radius and/or center discrepancy]
    \label{lem:time_diff}
    Let $\mathcal{B}_1$ and $\mathcal{B}_2$ be two Wasserstein balls of order $\pwass$, centered at arbitrary distributions $\mathbf{P}_1$ and $\mathbf{P}_2$, and with radii $\epsnom_1, \epsnom_2 > 0$.  For any  $x \in \mathcal{X}$,
\begin{align*}
        & \max_{\mathbf{Q} \in \mathcal{B}_1} \mathbf{E}_{\mathbf{Q}}[f(\unc,x)] -    \max_{\mathbf{Q} \in \mathcal{B}_2} \mathbf{E}_{\mathbf{Q}}[f(\unc,x)]  \leq L(\epsnom_1 - \epsnom_2 )_+ +  F_{\pwass} W_\pwass(\mathbf{P}_1,\mathbf{P}_2),
\end{align*}
where $F_{p} = L$ if $p=1$ and $F_{p} = \max_{x\in \mathcal{X}}(\int \|\nabla f(u,x)\|_*^qd\mathbf{Q})^{1/q}$ if $p>1$.
\end{lemma}
For $\pwass=1$ this outperforms the classic Kantorovich--Rubinstein bound~\citep{kantorovich_duality}, since then $F_p = L$ and $(\epsnom_1 - \epsnom_2)_+ \leq \epsnom_1 + \epsnom_2$; for $\pwass > 1$ the two are not directly comparable, as $F_p \leq L$ while $W_\pwass \geq W_1$.
This bound then motivates the following definition.
\begin{definition}[Worst-case expectation gap]
\label{def:gap}
    For two Wasserstein balls $\mathcal{B}_1$ and $\mathcal{B}_2$ of order $\pwass$, with centers $\mathbf{Q}_1, \mathbf{Q}_2 \in \Xi_\pwass(\supp)$ and radii $\epsnom_1, \epsnom_2 \geq 0$, define the (ordered) expectation gap
   $$\gap(\mathcal{B}_1,\mathcal{B}_2) = \min\left\{L(\epsnom_1 + \epsnom_2 +  W_1(\mathbf{Q}_1,\mathbf{Q}_2)),~  L(\epsnom_1 - \epsnom_2 )_+ +  F_{\pwass} W_\pwass(\mathbf{Q}_1,\mathbf{Q}_2) \right\},$$
   where a distribution $\mathbf{Q}$ in place of an ambiguity set denotes the degenerate set $\{\mathbf{Q}\}$ of radius $0$.
\end{definition}
The two clauses of $\gap$ come from Kantorovich duality and Lemma~\ref{lem:time_diff}, respectively.
Next, we also quantify the discrepancies induced by the clustering itself.
\begin{definition}[Clustering discrepancies]
\label{def:cluster_value}
    For any clustering of the dataset $\mathcal{D}_t$, with clusters $C^k_t$, cluster means $\bar{u}^k_t$, and clustered empirical distribution $\hat{\prob}^K_t$, define the clustering cost $c_{t,\pwass}$ and the piece-switching error $e_t$,
    $$c_{t,\pwass} = \left(\frac{1}{n_t}\sum_{k=1}^{K_t}\sum_{\hat{\unc} \in C^k_t}\left\| \hat{\unc} - \bar{u}^k_t\right\|_2^{\pwass}\right)^{1/{\pwass}}, \quad  e_t = \frac{1}{n_t}\sum_{k=1}^{K_t}\sum_{\hat{\unc} \in C^k_t}\max_{j\leq J}\, (z_t^{{jk}})^T(\bar{\unc}^k_t - \hat{\unc}),$$
    where $z^{jk}_t$ are dual variables from~\eqref{eq:robustopt_p_max}, interpretable as supergradients of the pieces $f_j$ in $u$ at the worst-case scenarios (for maximum-of-affine $f$, exactly the slopes $a_j(x^K_t)$).
    The error $e_t$ thus averages, at no extra computational cost, a first-order estimate of the worst per-piece change incurred by replacing each datapoint with its cluster mean; when $J=1$, deviations from the mean cancel within each cluster and $e_t = 0$, while the maximum prevents this cancellation for $J\geq2$.
\end{definition}
We are then ready to present a static (same time step) difference between the online and nominal~\gls{DRO} performance, under the following general assumptions.

\begin{assumption}
\label{ass:general_ass} Fix an order $\pwass$ and shrinking radius schedule $\epsnom_t$. For all $t \geq 0$:
    \begin{enumerate}
        \item The distribution $\hat{\prob}^K_t$ is constructed using cluster supports $\mathcal{S}_t$ that satisfy Assumption~\ref{ass:support}.
        \item The adaptive ambiguity set $\mathcal{P}^K_{t}$ has radius $\eps_{t} = \epsnom_{t}$.
    \end{enumerate}
\end{assumption}
Assumption~\ref{ass:general_ass} holds throughout the remainder of the paper.
\begin{lemma}[Static bounds] 
\label{thm:nominal_performance}
The value $H^K_t$ of the compressed \gls{DRO} problem~\eqref{eq:dro:time} satisfies
\ifpreprint
\begin{equation*}
\label{eq:relate}
\begin{gathered}
-\underline{\psi}_t  \le H_t^K - H_t \le \bar{\psi}_t + \delta,\\
\underline{\psi}_t = \min\left\{e_{t}, \gap(\mathcal{P}_t,\mathcal{P}^K_t)\right\}, \qquad \bar{\psi}_t = \min\left\{ L\epsnom_t +{\frac{M}{2}}(c_{t,2})^2,  \gap(\mathcal{P}_t,\mathcal{P}^K_t) \right\},
\end{gathered}
\end{equation*}
\else
\begin{equation*}
\label{eq:relate}
-\underline{\psi}_t  \le H_t^K - H_t \le \bar{\psi}_t + \delta, \quad \underline{\psi}_t = \min\left\{e_{t}, \gap(\mathcal{P}_t,\mathcal{P}^K_t)\right\}, ~~\bar{\psi}_t = \min\left\{ L\epsnom_t +{\frac{M}{2}}(c_{t,2})^2,  \gap(\mathcal{P}_t,\mathcal{P}^K_t) \right\},
\end{equation*}
\fi
where $H_t$ is the nominal~\gls{DRO} value and $\delta$ is the tolerance.  The bounds $\bar{\psi}_t$ and $\underline{\psi}_t$ have inner terms given in Definitions~\ref{def:gap} and~\ref{def:cluster_value}.
If smoothness (Assumption~\ref{ass:general_ass_smooth}) does not hold, $M=\infty$.  
\end{lemma}

The guarantees of nominal \gls{DRO} are thus inherited by our online solutions, with certificates adjusted by the clustering-based discrepancies $\underline{\psi}_t$ and $\bar{\psi}_t$, whose best clause can be chosen depending on the function $f$.
The clause involving $W_{1}(\hat{\prob}_t,\hat{\prob}^K_t)$, a Wasserstein-1 distance between two empirical distributions, is always computable via a linear program, and under smoothness, the $c_{t,2}$ term is applicable and easily computed from the given empirical distributions.
Furthermore, in some cases these bounds simplify: when $J=1$, \ie, $f$ consists of a single concave function, we have $\underline{\psi}_t = e_t = 0$. When instead $f$ is maximum-of-affine, $M = 0$, so $\bar{\psi}_t = \min\{L\epsnom_t,\gap(\mathcal{P}_t,\mathcal{P}^K_t)\}$. 

We use this result to formulate a deterministic bound $\tilde{R}(T,K)$ for the dynamic regret~\eqref{regret}.

\begin{theorem}[Deterministic regret]
\label{lem:regret}
    We have
    \begin{equation}
        \label{eq:regret_det}
    \begin{aligned}
    R(T,K) &\leq  \frac{1}{T}\sum_{t=0}^{T-1} \left(\underline{\psi}_{t+1} + \bar{\psi}_t + \gap(\mathcal{P}^K_{t+1},\mathcal{P}^K_t) \right) + \frac{1}{T} \gap(\mathcal{P}_0,\mathcal{P}_T) +\delta = \tilde{R}(T,K), \\
    \end{aligned}
    \end{equation} 
    with $\underline{\psi}_{t}$ and $\bar{\psi}_t$ as defined in Lemma~\ref{thm:nominal_performance}, and $\gap(\cdot)$ as defined in Definition~\ref{def:gap}.
    \end{theorem}
The bound has three sources of discrepancy: clustering, captured by $\underline{\psi}_{t+1} + \bar{\psi}_t$; the difference in time step, captured by the $\gap$ terms and the index offset in $\underline{\psi}_{t+1}$; and the optimality tolerance $\delta$.
Under clustering consistency (Lemma~\ref{lem:cluster_converge}), the clustering discrepancy converges to a fixed value and the time discrepancy to $0$, which we analyze next.
Specifically, the clustering-based values we previously defined have the following asymptotic limits.
\begin{lemma}[Asymptotic Wasserstein limit]
\label{thm:clustering_conv}
With $\prob^K_\star$ as in Lemma~\ref{lem:cluster_converge}, we have
$\mathbf{P}^\infty\{  w_K(\prob) \leq \lim_{t \rightarrow \infty }  W_{\pwass}(\hat{\prob}_t,\hat{\prob}^K_t) =W_{\pwass}(\prob,\prob^K_\star)\} = 1$,
where $ w_K(\prob) $ is the minimum quantization error of $\prob$.
Furthermore, $\lim_{K \rightarrow \infty} W_{\pwass}(\prob,\prob^K_\star) = 0$.
\end{lemma}
\begin{lemma}[Asymptotic clustering-discrepancy limits]
\label{thm:asy_cluster}
Let $\kappa(u)$ denote the index of the fixed cluster support $S^{\kappa(u)}_\tau$ containing $u$, which is $\prob$-almost surely unique by Assumption~\ref{ass:support}.
With $\prob^K_\star$ as in Lemma~\ref{lem:cluster_converge}, we have
$\prob^{\infty}\{\lim_{t\to\infty} c_{t,\pwass} = c_{\star,\pwass}= (\mathbf{E}_{\prob}[\|u - \bar{u}^{\kappa(u)}_\star\|^{\pwass}_2 ])^{1/{\pwass}} \} = 1$, with $\bar{u}^k_\star$ the limits of the cluster means (and are therefore atoms of $\prob^K_\star$). 
Furthermore, if there exist accumulation points of $z_t^{jk}$ for all $j\leq J,k\leq K$, then
$\prob^{\infty}\{\limsup_{t\to\infty}e_t = e_\star = \mathbf{E}_{\prob}[\max_{j\leq J}(z^{j\kappa(u)}_\star)^T(\bar{u}^{\kappa(u)}_\star - u) ]\} = 1$,
where $z^{jk}_\star$ are the accumulation points that achieve the limit superior.
\end{lemma}
We are now ready to present our asymptotic regret, which shows that as $T \rightarrow \infty$, the time-discrepancy disappears, and the clustering discrepancy is bounded by the above limit terms.
\begin{theorem}[Asymptotic regret]
\label{cor:regret}
Let $\epsnom_t$ and $\beta_t$ be chosen according to Theorem~\ref{thm:calc_eps}. 
With probability 1,
   $ \lim_{T\to\infty} R(T,K) \leq O(\max\{\underline{\psi}_\star,\bar{\psi}_\star\})+\delta,$ 
   where 
   $\underline{\psi}_\star = \min \left\{
   \limsup_{t \rightarrow \infty } {e_t}, \gap(\prob,\prob^K_\star) \right\},~
        \bar{\psi}_\star =
\min \left\{ \frac{M}{2}(c_{\star,2})^2,  \gap(\prob,\prob^K_\star)\right\},~ 
\gap(\prob,\prob^K_\star) = \min\left\{ LW_{1}(\prob,\prob^K_\star), F_{p}W_{\pwass}(\prob,\prob^K_\star)\right\}$.
   Additionally, as $K \to \infty$,~$\lim_{K\to\infty}\lim_{T\to\infty}R(T,K) = 0.$
\end{theorem}

For maximum-of-affine $f$, the following corollary additionally simplifies some of the clauses.
\begin{corollary}
    \label{cor:regret_opt}
    Under the assumptions of Theorem~\ref{cor:regret}, if the function $f(u,x)$ is maximum-of-affine in both $u$ and $x$, with $\supp = \reals^d$, we have $\limsup_{t \to \infty} e_{t} = e_\star$
where the latter is defined in Lemma~\ref{thm:asy_cluster}, and $\bar{\psi}_\star = 0$.
Furthermore, if the clustering-induced coupling is optimal with respect to $\pwass=1$, \ie,
    $ c_{t,1} = W_{1}(\hat{\prob}_t,\hat{\prob}^K_t), $
   then almost surely with respect to~$\prob^{\infty}$,
   $\underline{\psi}_\star \leq {e_\star}\leq
       LW_{1}(\prob,\prob^K_\star).$
\end{corollary}

Overall, the converged values always depend on the number of clusters $K$ and the clustering algorithm, which control the quality of the clustered distribution and its limit $\prob^K_\star$, and hence the trade-off between efficiency and optimality.

So far, we have obtained deterministic regret bounds. We can, however, also establish probabilistic bounds, which gives the {\it rate} at which the regret bound converges to the asymptotic limit above. This requires an additional lemma on the clustering convergence rate. 

\begin{lemma}[Clustering convergence rate]
\label{lem:conditional_Conv}
Suppose the sequence $\epsnom_t$ is computed using Theorem~\ref{thm:calc_eps} such that $\lim_{t \rightarrow \infty} \epsnom_t = 0$ and $\beta_t \in (0,1)$ satisfies $\sum_{t=1}^\infty {\beta_t < \infty}$. Then, for $t \geq \tau$, the distance $W_{\pwass}(\hat{\prob}^K_t,\prob^K_\star)$ converges to $0$ with rate $O(\epsnom_t + (\log(\beta_t^{-1})/t)^{1/(2\pwass)})$. For $\pwass = 1$, this rate is $O(\epsnom_t)$.
\end{lemma}

\begin{theorem}[Dynamic regret bound]
\label{thm:regret}
Let $\epsnom_t$ and $\beta_t$ be chosen according to Theorem~\ref{thm:calc_eps}. For finite $T \gg \tau$, with probability $1-\bar{\beta}$, the dynamic regret is bounded as 
   $R(T,K) \leq O\left(\left(\log(\beta_T^{-1})/T\right)^{\min(\pwass/d,1/2)}\right) + 
O\left(W_1(\prob,\prob^K_\star)\right) + \delta,$
   where $\beta_T$ decays at least sublinearly in $T$, and
    $1-\bar{\beta} \geq \Pi_{t=1}^T(1-\beta_t)^4$.
\end{theorem}

This provides a sublinear convergence rate for the dynamic regret, matching the convergence rate of the underlying empirical distributions, for which no sharper rates exist in the literature.

\section{Subgradient method}
\label{sec:subgrad_sec}
When the~\gls{DRO} problem is too expensive to solve in full, we can instead obtain an approximate solution.
In this section, for general $p$, we derive fast decompositions of the worst-case expectation, $\max_{\mathbf{Q} \in \mathcal{P}^K}\mathbf{E}_{\mathbf{Q}}[f(\unc,x^K_{t-1})]$, under various supports $S$ and ground norms $\|\cdot\|$ of the dynamic ambiguity set $\mathcal{P}^K$.
The results differ depending on the finiteness of $p$, so for $p=\infty$, we instead denote the ambiguity set $\mathcal{P}^K_\infty$.
We focus on $f$ that is maximum-of-affine in $u$, given as $$f(u,x) = \max_{j \leq J} (a_j(x)^T u + b_j(x)),$$ where $a_j(x) \in \reals^d $ and $b_j(x)\in \reals$ are convex functions of $x$. For simplicity, throughout this section we denote $\bar{x} = x^K_{t-1}$, drop subscripts $t$, and let $\epsnom = \eps$.

\subsection{Time complexity}
\label{sec:general_structure}
We reduce the worst-case expectation by decomposing it into per-centroid and per-piece oracles, and construct the subgradient from the individual solutions. We begin with the following definitions, which hold by strong duality~\citep{mohajerin2018data, wang2022mean}.
\begin{definition}
    \label{def:oracles}
Fix $\bar{x}$. For a centroid $u\in S$ and piece $j\leq J$, denote the per-piece oracles:
$$\varphi_\lambda^j(u) \define \sup_{\zeta\in S}\big(a_j(\bar x)^T\zeta+b_j(\bar x)-\lambda\|\zeta-u\|^p\big),\ \ \lambda\geq0,
\qquad
\varphi_\infty^j(u) \define \sup_{\zeta\in S,~\|\zeta-u\|\leq\epsnom}\big(a_j(\bar x)^T\zeta+b_j(\bar x)\big),$$
with $p \geq 1$ (left) and $p=\infty$ (right). The respective per-centroid oracles are $\varphi_\lambda(u)\define\max_{j\leq J}\varphi_\lambda^j(u)$ and $\varphi_\infty(u)\define\max_{j\leq J}\varphi_\infty^j(u)$, and the worst-case expectations are given by
$$\max_{\mathbf{Q} \in \mathcal{P}^K}\mathbf{E}_{\mathbf{Q}}[f(u,\bar{x})] = \inf_{\lambda\geq0}h(\bar{x},\lambda), \qquad \max_{\mathbf{Q} \in \mathcal{P}^K_\infty}\mathbf{E}_{\mathbf{Q}}[f(u,\bar{x})] = \sum_{k=1}^{K}\theta^k\varphi_\infty(\bar u^k),$$
where $h(\bar{x},\lambda) \define \lambda{\epsnom}^p + \sum_{k=1}^{K}\theta^k\varphi_\lambda(\bar u^k)$.
\end{definition}

At optimality, the subgradient takes the following form. 
\begin{proposition}
    \label{thm:gen_subgrad}
Let $\lambda^\star(\bar{x})\in\argmin_{\lambda\geq0}h(\bar{x},\lambda)$, and let $\zeta^{k}$ be a maximizer of cluster $k$'s oracle at $\bar u^k$ ($\varphi_{\lambda^\star}(\bar u^k)$ or $\varphi_\infty(\bar u^k)$).
Except when $p=1$ with $S = \reals^d$, where Corollary~\ref{cor:p1_closed} gives an explicit subgradient, a subgradient of the worst-case expectation at $\bar x$ is
\begin{equation}
\label{eq:general_subgrad}
g \in \sum_{k=1}^{K}\theta^k\,\partial_x f(\zeta^{k},\bar{x}).
\end{equation}
\end{proposition}
As presented, this subgradient may not be unique, but is exact. However, in order to speed up the computation, for general bounded support $S$, we may instead solve the worst-case expectation to a value-tolerance $\delta_{\text{approx}}$, which also subjects the subgradient to a tolerance $\delta_{\text{grad}}$, defined below.

\begin{lemma}[Subgradient tolerance]
\label{thm:combined_tol}
Suppose $S$ is bounded, with $R_S \define \sup_{\zeta,\zeta'\in S}\|\zeta-\zeta'\| < \infty$. For finite $p\geq1$, suppose $\hat\lambda\geq0$ satisfies $|\hat\lambda-\lambda^\star(\bar x)|\leq\sigma$, and that for each $k\leq K$, the inner problem at $\hat\lambda$ is solved to value-tolerance $\Gamma$, at solution $\hat\zeta^k \in S$. Then, the worst-case expectation is solved to value-tolerance $\delta_{\text{approx}} = \max\{\Gamma,R_S^p\sigma\}$.
Furthermore, for any $x \in \mathcal{X}$, any subgradient
$\hat g \in \sum_{k=1}^{K}\theta^k\,\partial_x f(\hat\zeta^k,\bar x)$ satisfies
\begin{equation}
\label{eq:approx_subgrad_ineq}
\max_{\mathbf{Q}\in\mathcal{P}^K}\mathbf{E}_{\mathbf{Q}}[f(u,x)] \geq \max_{\mathbf{Q}\in\mathcal{P}^K}\mathbf{E}_{\mathbf{Q}}[f(u,\bar x)] + \hat g^T(x-\bar x) - \delta_{\text{grad}},
\end{equation}
where $\delta_{\text{grad}} = \Gamma + R_S^p\sigma + R_S^p\big|\lambda^\star(\bar x)-\lambda^\star(x)\big|$.
If $p=\infty$, $\delta_{\text{approx}} = \delta_{\text{grad}} = \Gamma$.
\end{lemma}
The subgradient tolerance is required for quantifying the additional regret we incur with inexact subgradient solves. We now derive the computational complexities of our subroutines for obtaining $g$ (when the solution is closed-form) or $\hat{g}$, which are notably faster than the cost $O((KJd)^{1.5})$ (or $O((KJ(m+d))^{1.5})$ for polyhedral $S$ with $m$ facets) of the original formulation.

\begin{proposition}[Time complexity]
\label{thm:finitep}
If $S=\reals^d$, the oracles admit closed form solutions at cost $T_{\rm inner}=O(1)$. 
If instead $S\subsetneq\reals^d$, the oracles follow the subroutines of Table~\ref{tab:inner_cost} at cost $T_{\rm inner}(S,\|\cdot\|,p)$.
For finite $p$, except when $p=1$ with $S = \reals^d$, $\hat{\lambda}$ is found through bisection in $T_{\rm outer}=O(\log(1/\sigma))$ iterations, with tolerance $\sigma$; for $p=\infty$, no outer search is needed since no shared $\lambda$ exist, and $T_{\rm outer}=1$.
The subgradient is computable in total time $O(T_{\rm outer}KJ\,T_{\rm inner}(S,\|\cdot\|,p))$, or $O(T_{\rm outer}\,T_{\rm inner}(S,\|\cdot\|,p))$ if the $KJ$ inner solves are parallelized.
\end{proposition}
\begin{corollary}[Closed form at $p=1$, unbounded $S$]
\label{cor:p1_closed}
At $p=1$ and $S=\reals^d$, let $j^k\in\argmax_{j\leq J}\big(a_j(\bar x)^T\bar u^k+b_j(\bar x)\big)$, and $j^\star\in\argmax_{j\leq J}\|a_j(\bar{x})\|_*$.
The subgradient is
$g \in \sum_{k=1}^{K}\theta^k\,\partial_x\big(a_{j^k}(\bar x)^T\bar u^k+b_{j^k}(\bar x)\big)
+ \epsnom\,\partial_x\|a_{j^\star}(\bar{x})\|_*.$
This requires only $K+1$ maximizations over the $J$ pieces and no bisection over $\lambda$, since $\lambda^\star(\bar{x}) = \|a_j(\bar{x})\|_*$.
\end{corollary}

\begin{corollary}[Closed form at $p=\infty$, unbounded $S$]
\label{cor:pinf_closed}
At $p=\infty$ and $S=\reals^d$, let $j^k\in\argmax_{j\leq J}\big(a_j(\bar x)^T\bar u^k+b_j(\bar x)+\epsnom\|a_j(\bar x)\|_*\big)$.
The subgradient is
$g \in \sum_{k=1}^{K}\theta^k\Big(\partial_x\big(a_{j^k}(\bar x)^T\bar u^k+b_{j^k}(\bar x)\big)+\epsnom\,\partial_x\|a_{j^k}(\bar{x})\|_*\Big).$
This requires only $K$ maximizations over the $J$ pieces.
\end{corollary}

Table~\ref{tab:inner_cost} summarizes the inner-problem subroutines and costs, Example~\ref{ex:box_l2_p1} gives one derivation, and Appendix~\ref{sec:inner_derivations} gives the rest.
We note that translating an inner bisection tolerance $\sigma_{\text{in}}$ into the value-tolerance $\Gamma$ of Lemma~\ref{thm:combined_tol} requires an additional, problem-specific Lipschitz argument, similar to the one for $\hat{\lambda}$, which we omit here for brevity.

\begin{table}[h]
\centering
\small
\caption{Inner problem subroutines and costs (Proposition~\ref{thm:finitep}) under various supports and ground norms, for finite $p\geq1$ and $p=\infty$; slashed entries separate the parenthesized special case from the general case. Here $d$ is the dimension of $u$, $m$ the number of facets of a polyhedral support, and $\sigma_{\text{in}}$ the inner bisection tolerance.}
\label{tab:inner_cost}
\begin{adjustbox}{max width=\textwidth}
\begin{tabular}{llllll}
\toprule
& & \multicolumn{2}{l}{Finite $p\geq1$} & \multicolumn{2}{l}{$p=\infty$}\\
\cmidrule(lr){3-4}\cmidrule(lr){5-6}
Support $S$ & Norm & Subroutine & Cost & Subroutine & Cost\\
\midrule
Unbounded, $\reals^d$ & any & Closed form & $O(1)$ & Closed form & $O(1)$\\
Box & $\|\cdot\|_1$ & Threshold ($p=1$) / sort + scan & $O(d)$ / $O(d\log d)$ & Sort + fill & $O(d\log d)$\\
Box & $\|\cdot\|_\infty$ & Sort + scan & $O(d\log d)$ & Clip & $O(d)$\\
Box & $\|\cdot\|_2$ & Clip ($p=2$) / bisection & $O(d)$ / $O(d\log(1/\sigma_{\text{in}}))$ & Bisection & $O(d\log(1/\sigma_{\text{in}}))$\\
Polyhedral & $\|\cdot\|_1$ & LP ($p=1$) / power cone & $O(d(m+d)^{1.5})$ & LP & $O(d(m+d)^{1.5})$\\
Polyhedral & $\|\cdot\|_2$ & SOCP ($p\in\{1,2\}$) / power cone & $O((m+d)^{1.5})$ & SOCP & $O((m+d)^{1.5})$\\
Norm ball & $\|\cdot\|_2$ & Bisection & $O(d\log(1/\sigma_{\text{in}}))$ & Bisection & $O(d\log(1/\sigma_{\text{in}}))$\\
\bottomrule
\end{tabular}
\end{adjustbox}
\end{table}

\begin{example}[Box support, $\ell_2$ norm, $p=1$]
\label{ex:box_l2_p1}
Fix a box $S=[l,h]$, a point $u\in S$, and $\lambda\geq0$. Simplify $\beta\define a_j(\bar x)$, and consider
$\varphi_\lambda^j(u)=\sup_{\zeta\in S}\big(\beta^T\zeta-\lambda\|\zeta-u\|_2\big).$

\textit{Case $\|\beta\|_2\leq\lambda$.} By the Cauchy--Schwarz inequality, $\zeta^\star=u$ is optimal. Cost $O(1)$.

\textit{Case $\|\beta\|_2>\lambda$.} 
Define $\mu = \lambda/\|\zeta^\star-u\|_2$.
The gradient of the penalty $\lambda\|\zeta-u\|_2$ w.r.t. $\zeta$ is $\mu(\zeta^\star-u)$ at optimality. Then, solving the KKT system w.r.t. the box support, we have $\zeta_i^\star(\mu)=\mathrm{clip}(u_i+\beta_i/\mu,\,l_i,\,h_i)$, a function of $\mu$.
To solve for $\mu$, note that it must satisfy the condition $\mu\|\zeta^\star(\mu)-u\|_2 = \lambda$. Since $\|\zeta^\star(\mu)-u\|_2$ is continuous and monotonically increasing from $0$ to $\|\beta\|_2$, the function $\mu\|\zeta^\star(\mu)-u\|_2 - \lambda$ has a unique root $\mu^\star$, which can be found by bisection. With bisection tolerance $\sigma_{\text{in}}$, the cost is $O(d\log(1/\sigma_{\text{in}}))$.
\end{example}

\subsection{Adjustments to the dynamic regret}
\label{sec:regret_sub}
Throughout Section~\ref{sec:regret}, $x^K_t$ is the minimizer of the compressed problem, solved to tolerance $\delta$. We now instead bound the regret of the {\it subgradient} solutions $x^K_t$, with the worst-case expectation solved via the above subroutines.
For this, we require that $f$ is Lipschitz in $x$ w.r.t.\ the $\ell_2$ norm.
\begin{assumption}[Lipschitzness in $x$]
\label{ass:subgrad_lip}
For every $\unc \in \supp$ and $j \leq J$, the function $f_j$ is $L_{x,j}$-Lipschitz on $\mathcal{X}$ w.r.t.\ the $\ell_2$-norm. Let $L_x = \max_{j\leq J} L_{x,j}$.
\end{assumption}
Let $R_{\text{approx}}(T,K)$ denote the dynamic regret~\eqref{regret} evaluated at the subgradient iterates $x^K_t$.
Compared to the regret of Theorem~\ref{lem:regret}, we incur three additional discrepancies: (i) the standard convergence rate of the subgradient method on a fixed objective; (ii) the drift of the objective itself, as $\hat{\prob}^K_t$ shifts and the radius $\eps_t$ shrinks; and (iii) the subgradient tolerance from the subroutines of Section~\ref{sec:general_structure}.
\begin{theorem}[Subgradient dynamic regret]
\label{thm:regret_single_resolve}
Suppose Assumption~\ref{ass:subgrad_lip} holds. Let $\mathcal{X}$ be convex, and $R_{\mathcal{X}} = \mathrm{diam}_2(\mathcal{X})$. 
Suppose $x^K_t$ is obtained through the subgradient procedure of Algorithm~\ref{alg:hl}, with step size $\eta_t = R_{\mathcal{X}}/(L_x\sqrt{t})$.
Then, with  $\tilde{R}(T,K)$ the full-solve regret bound~\eqref{eq:regret_det},
\begin{equation}
\label{eq:regret_single_resolve}
R_{\text{approx}}(T,K)  \leq \tilde{R}(T,K) + \frac{3R_{\mathcal{X}}L_x}{2\sqrt{T}} + \frac{1}{T}\Big(\gap(\mathcal{P}^K_T,\mathcal{P}^K_0)+ 2\sum_{t=1}^T \gap(\mathcal{P}^K_0,\mathcal{P}^K_t)\Big) +\frac1T\sum_{t=1}^T\delta_{\text{grad},t},
\end{equation}
where each $\gap$ term (Definition~\ref{def:gap}) measures the drift of the worst-case expectation between the ambiguity sets at two timesteps,
and $\delta_{\text{grad},t}$ is the subgradient tolerance of Lemma~\ref{thm:combined_tol}, evaluated at $\bar x=x^K_{t-1}$ and $x=x^\star$, a minimizer of the compressed problem at time $0$.
\end{theorem}
The first additional term, the convergence rate of the subgradient method, vanishes as $T$ grows, while the other two additional terms do not.
However, if the subgradient is exact, the last term disappears.
Furthermore, if the problem is periodically re-solved in full, each drift term $\gap(\mathcal{P}^K_0,\mathcal{P}^K_t)$ is replaced by $\gap(\mathcal{P}^K_{\ell(t)},\mathcal{P}^K_t)$, where $\ell(t) \leq t$ is the latest resolve time.
Then, as the clustered distribution and radius converge, the resolve-relative drifts vanish, and the entire second additional term also converges to 0. We would then recover the direct-solve regret.

\section{Online algorithms for data-compression}
\label{sec:algos}
In Section~\ref{sec:shared_cluster}, we gave clustering requirements, and in this section, we provide two online clustering algorithms satisfying Assumption~\ref{ass:support}.
These algorithms maintain a reasonably accurate approximation of the optimal discretization while minimizing computational effort.

\subsection{Approximate $k$-centers clustering with warm starts (Reclustering)}
\label{sec:recluster}
To maintain a consistently good approximation of the data, we recompute an approximate $k$-centers clustering at each time step $t$, until the stopping threshold $\tau$, warm-starting each clustering with the previous centers $A_{t-1}$.
The support of cluster $k$ is the Voronoi region around its center $a_t^k \in A_t$.
After time $\tau$, the centers are fixed at $A_\tau$: the datapoint $\hat{u}$ observed at time $t$ joins the cluster $k'=\argmin_{k\leq K} \|\hat{u} - a^k_\tau\|$, the cluster mean is updated to $\bar{u}^{k'}_t = (n^{k'}_{t-1}\bar{u}^{k'}_{t-1} + \hat{u})/n^{k'}_t$, and the weights rescale to $\theta^k_{t} = (n^k_{t-1} + \ones\{k = k'\})/n_{t}$ for $k=1,\dots, K$.
We approximate the $k$-centers problem with $k$-means~\citep{kmeans}, whose iterations cost $O(n_t K d)$ time each; in practice, warm-starting greatly reduces the number of iterations.
This method is not memory efficient, however, as storing the datapoints to recompute the clusters requires $O((K + n_t)d)$ memory.

\subsection{Online clustering (OnlineClustering)}
\label{sec:clustream}
For an algorithm that is also memory conservative, we propose an online updating procedure that takes inspiration from CluStream~\citep{clustream}.
We maintain two sets of clusters: $Q \geq K$ {\it microclusters}, and $K$ {\it macroclusters}.
Each arriving datapoint either joins the microcluster whose support contains it, if within a distance threshold of its center, or starts a new microcluster; in the latter case, the two closest microclusters are merged to maintain the budget $Q$.
Every 50 timesteps, $K$ macroclusters are obtained by a warm-started approximate $k$-centers clustering of the microcluster centers.
Appendix~\ref{app:onlinecluster} details the procedure and shows that it satisfies the clustering requirements of Section~\ref{sec:shared_cluster}.
Note that the clustering iterations cost $O(QKd)$, ordinary microcluster updates cost only $O(Qd)$, and the memory requirement is $O((K+Q)d)$, so, unlike Reclustering, neither time nor memory grows with $n_t$.

\section{Numerical examples}
\label{sec:exampless}
We now illustrate the computational performance and robustness of the proposed method on two numerical examples: portfolio optimization and sparse \gls{SVM}.
All the code to reproduce our experiments is available, in Python, at 
\begin{center}
{\tt \url{https://github.com/stellatogrp/online_mro}}.
\end{center}
We run the experiments on the Princeton Research Computing facility, on jobs with up to 32 parallel cores.
We solve all mixed-integer optimization problems with the MOSEK~\citep{mosek} optimizer with default settings, and all continuous problems with Clarabel.

We compare the approaches given in Table~\ref{tab:baselines} over a streaming training dataset of $T = 5000$ time steps.
The exact methods re-solve~\eqref{eq:robustopt_p_max} every 100 (portfolio) or 200 (\gls{SVM}) time steps and at $t \in T_{\text{full}}=\{1, 5, 25, 50, 100\}$.
The subgradient methods take one projected step per time step with step size $\eta_t = 0.1/\sqrt{t}$, re-anchored by a full solve every 200 time steps and at $t \in T_{\text{full}}$, warm-starting from the subgradient iterate.
As an oracle for the true optimal value, we solve the \gls{SAA} problem on $N = 20{,}000$ samples, shown as $H^\star$ in the figures.

We report six metrics, averaged over 30 repetitions of each experiment, with the 25--75\textsuperscript{th} percentiles shown as shaded regions: (i-ii) the in-sample and out-of-sample objective values, the latter being the empirical expectation of $f$ over the testing data; (iii) the empirical confidence $1 - \hat{\beta}_t$, \ie, the fraction of repetitions in which the in-sample value upper bounds the out-of-sample value; (iv) the per-iteration times of the exact methods, including both clustering and solving times for the data-compressed methods; (v) the cumulative times of the subgradient-based methods; and (vi) the dynamic regret~\eqref{regret} of online clustering, with its theoretical upper bound~\eqref{eq:regret_det}.

\paragraph{Hyper-parameter selection}
If the initial dataset is large enough, the number of clusters $K$ can be chosen a priori, using the elbow method on the clustering cost $(c_{0,2})^2$~\citep{elbow,wang2022mean}; otherwise, it can be initialized small and adapted before the freeze time $\tau$.
We use $K = 15$ in both experiments.
For the ambiguity set radius, each \gls{DRO} method sweeps a finite grid of initial radii, decayed polynomially over time.
For each method independently, we select the sequence whose empirical confidence over the last 10\% of time steps reaches at least $0.6$ and whose mean out-of-sample value over that window is smallest; if no candidate reaches the threshold, we select the sequence with the highest confidence.

\begin{table}[h]
  \centering
  \caption{The compared approaches.}
  \label{tab:baselines}
  \begin{tabular}{p{0.32\textwidth}p{0.6\textwidth}}
    \toprule
    Method & Description\\
    \midrule
    Online clustering (ours) & Solve~\eqref{eq:robustopt_p_max} w.r.t.\ $\hat{\prob}^K_t$, updated by our online clustering algorithm (Section~\ref{sec:clustream}).\\
    Subgradient online clustering (ours) & The above method, with subgradient solves (Section~\ref{sec:subgrad_sec}).\\
    Wasserstein \gls{DRO} & Solve~\eqref{eq:robustopt_p_max} with the ambiguity set $\mathcal{P}_{t}$.\\
    Subgradient Wasserstein \gls{DRO} (ours) & Wasserstein \gls{DRO} with subgradient solves (Section~\ref{sec:subgrad_sec}).\\
    Sample average approx. (\gls{SAA}) & Solve the stochastic optimization (SO) problem w.r.t. $\hat{\prob}_{t}$.\\
    Clustered \gls{SAA} &  Solve the \gls{SO} problem w.r.t. $\hat{\prob}^K_{t}$ obtained via $k$-means.\\
    \bottomrule
  \end{tabular}
  \vspace{-1em}
\end{table}

\subsection{Portfolio optimization}
\label{sub:sparse_portfolio_optimization}
We consider a market that forbids short-selling and has $d$ assets~\citep{wang2022mean,kuhn2019wasserstein}.
The uncertain parameters are the daily returns, given by $u \in \reals^d$.
The decisions ${x\in \reals^n}$ represent the percentage weights (of the total capital) invested in the assets.
The distribution $\prob$ is unknown, but we have access to a streaming dataset $\mathcal{D}_t$, updated each day with a new returns vector. 
Our objective is to minimize the $\CVaR$ with respect to $x$,
$$\text{minimize}~\CVaR(- u^Tx,\omega), ~~\text{subject to}~\ones^T x = 1,~ x \ge 0,$$
which represents the average of the $\omega$ largest portfolio losses that occur. In other words, the $\CVaR$ term seeks to ensure that the expected magnitude of portfolio losses, when they occur, is low.
The objective can be written as the expectation of the maximum-of-affine functions~\citep{cvaropt}, \ie, $\Expect_{\prob}\left[\gamma +(1/\omega) \max\{-u^Tx -\gamma, 0\}\right]$. 
We solve the \gls{DRO} problem with $\pwass=1$, $\|\cdot\|$ the $\ell_2$-norm, and a box uncertainty set $S = [l, h]$.

\paragraph{Problem setup}
We generate synthetic returns using the factor model $u_j = \mu_j + \sum_{r=1}^{3} F_{r}\,\beta_{j,r} + \epsilon_{j}$, for stocks $j = 1,\dots,d$. We set $\mu_j \sim U(-0.01,0.04)$, $F_r\sim \mathcal{N}(0,0.03^2)$, $\beta_{j,r} \sim U(0.3,1.2)$, and noise $\epsilon_j \sim \mathcal{N}(0,\sigma_j^2)$, where $\sigma_j \sim U(0.02,0.06)$. 
For each stock, we set the lower and upper bounds $l_j$ and $h_j$ as the 30\% and 70\% quantiles of a separate historical window of 5000 observations from the same model, mimicking price-limit bands inferred from history.
We let $\omega = 20\%$ and $d = 300$ stocks, initialize with a dataset of size $n_0 = 5$, and use $Q = 500$ microclusters.
The out-of-sample values are computed on a test set of 1000 datapoints.
The selected radii decay as $\epsnom_t \sim t^{-1/20}$.

\begin{figure}[t!]%
	\centering
      \includegraphics[width = \textwidth]{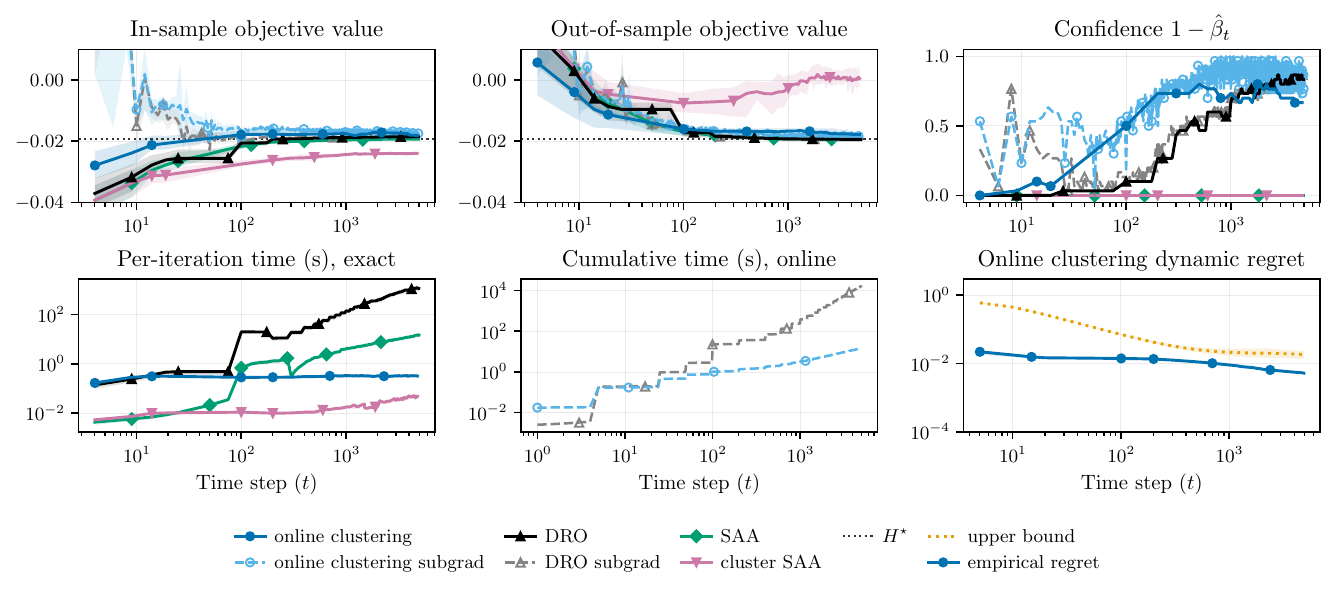}
	\caption{Portfolio optimization, $K=15$. Top: in-sample/out-of-sample objectives and confidence $1-\hat{\beta}_t$, with the oracle $H^\star$ dotted. Bottom: per-iteration times of the exact methods, cumulative times of the subgradient methods, and the dynamic regret of online clustering with its upper bound~\eqref{eq:regret_det}.}%
	\label{fig:port_combined}%
\end{figure}
\paragraph{Results}
Figure~\ref{fig:port_combined} shows the results.
In and out-of-sample, online clustering tracks nominal \gls{DRO} closely, with both approaching the oracle value $H^\star$. The compression introduces minimal visible loss.
The empirical confidence of the \gls{DRO}-based methods grows towards $0.8$, while \gls{SAA} and clustered \gls{SAA} never certify their solutions: their confidence remains near $0$ for all $t$, and clustered \gls{SAA}, a pure scenario reduction via $k$-means, is also unreliable out-of-sample. 
The per-iteration time of online clustering stays constant at roughly $0.3$ seconds, while nominal \gls{DRO} grows beyond $10^3$ seconds per solve and \gls{SAA} beyond $10$ seconds.
Among the subgradient variants, which apply the subroutines of Section~\ref{sec:subgrad_sec}, the compressed version accumulates orders of magnitude less time than the nominal one, whose per-step cost scales with $t$.
In this example, online clustering with exact or subgradient solves are the recommended methods: they nearly match the nominal \gls{DRO} values at a constant, small per-iteration cost.
The dynamic regret of online clustering decreases over time and remains below the theoretical upper bound, which converges to a $K$-dependent value, as given by Theorem~\ref{cor:regret}.

\subsection{Sparse SVM}
\label{sub:sparse_svm}
We consider a cardinality-constrained (best-subset) support vector machine for binary classification.
The uncertain parameter is a labeled example $\unc = (v,y)$, with covariate $v \in \reals^d$ and label $y \in \{-1,1\}$, and the decision $x \in \reals^d$ is a linear classifier restricted to at most $r$ nonzero coordinates.
The distribution $\prob$ is unknown, but we have access to a streaming dataset $\mathcal{D}_t$ of labeled examples.
Our objective is to minimize the worst-case hinge loss,
$$ \text{minimize}~\Expect_{\prob}\left[(1-y\,x^Tv)_+\right],~~
	\text{subject to}~\|x\|_0 \le r,$$
where the ambiguity set perturbs the covariate $v$ only, leaving the label $y$ fixed.
The hinge loss $(1-y\,x^Tv)_+$ is maximum-of-affine in $v$.
We use $\pwass=1$, the $\ell_1$ norm, and $S = \reals^d$; the cardinality constraint makes the exact problem a \gls{MILP}.

\paragraph{Problem setup}
We use the \texttt{ijcnn1} dataset from the LIBSVM repository, derived from the IJCNN 2001 neural network competition, which has $d=22$ real-valued covariates and $n=49{,}990$ labeled examples.
In each repetition, we draw a random train/test split, streaming $T = 5000$ training points and evaluating out-of-sample on a test set of roughly $8300$ points.
We set the cardinality budget $r=10$, initialize with a dataset of size $n_0 = 5$, and use $Q = 500$ microclusters. We make the clusters label-aware, such that datapoints are added only to same-label clusters. 
Since $\mathcal{X}$ is non-convex, the subgradient-based methods project each step onto the cardinality constraint via iterative hard thresholding~(IHT), with the periodic full solves acting as re-anchors.

\begin{figure}[t!]%
	\centering
      \includegraphics[width = \textwidth]{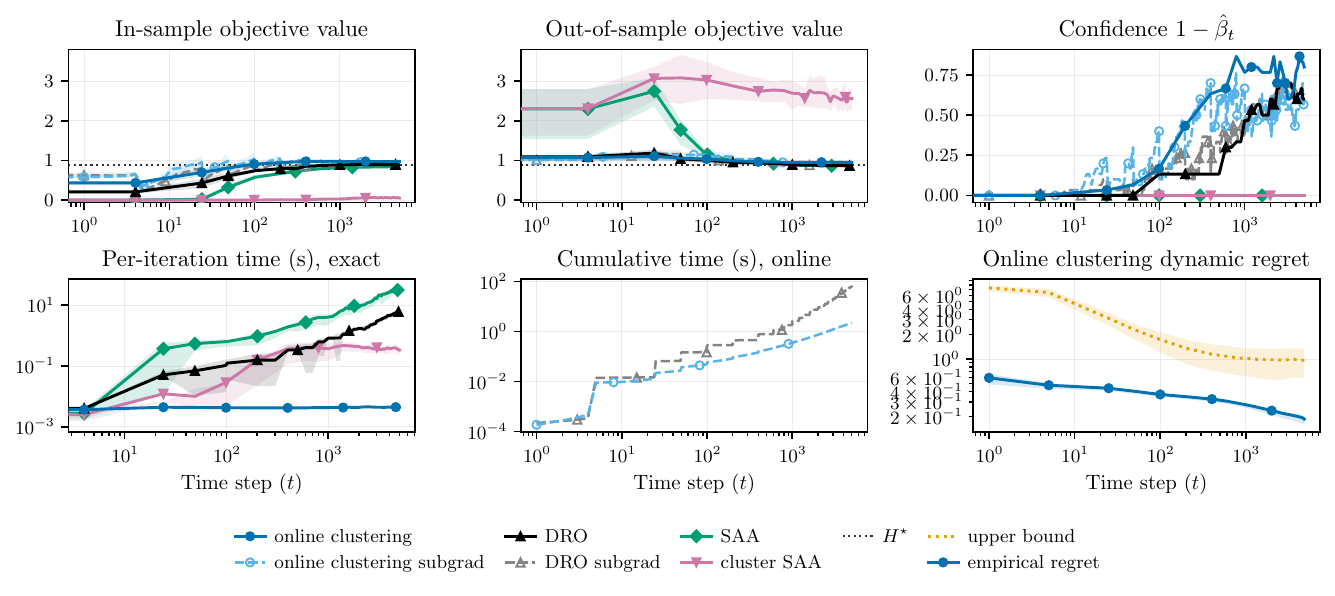}
	\caption{Sparse \gls{SVM}, $K=15$. Top: in-sample/out-of-sample objectives and confidence $1-\hat{\beta}_t$, with the oracle $H^\star$ dotted. Bottom: per-iteration times of the exact methods, cumulative times of the subgradient methods, and the dynamic regret of online clustering with its upper bound~\eqref{eq:regret_det}.}%
	\label{fig:svm_combined}%
\end{figure}
\paragraph{Results}
Figure~\ref{fig:svm_combined} shows the results.
The data-compressed approaches incur a small in-sample gap relative to nominal \gls{DRO}, but match it out-of-sample, approaching the oracle value $H^\star$. \gls{SAA} joins them only after a long time, and clustered \gls{SAA} remains roughly three times worse out-of-sample.
As in the portfolio example, only the \gls{DRO}-based methods certify their solutions, with confidence growing towards $0.8$.
\gls{SAA} is the most expensive exact method here, even beyond nominal \gls{DRO}; the \gls{DRO} problems act as {\it regularized} versions of the \gls{SAA} problem, which stabilizes the branch-and-bound search.
The per-iteration time of online clustering remains constant at roughly $4$ milliseconds, since the \gls{MILP} size is fixed by the $K=15$ scenarios.
Among the subgradient variants, the compressed version accumulates over an order of magnitude less time than the nominal one, and is therefore preferred. 
The dynamic regret of online clustering again decreases and stays below its upper bound.

\section{Conclusions}
\label{sec:conclusions}
In this paper, we introduced an online data-compression framework for solving Wasserstein~\gls{DRO} problems with streaming data.
By clustering the data online, our method keeps the problem size bounded as data accumulates, while inheriting the out-of-sample guarantees of nominal~\gls{DRO} with certificates adjusted by the compression discrepancies.
We quantified the price of compression through static performance bounds and dynamic regret bounds that converge sublinearly to a clustering-dependent gap, and extended the analysis to fast subgradient updates built on specialized worst-case-expectation subroutines.
Numerical experiments on portfolio optimization and sparse~\gls{SVM} problems, including mixed-integer formulations, show over an order of magnitude reduction in cumulative computation time with minimal loss in solution quality.

\ifpreprint
\section*{Acknowledgments}
\mythanks{}
\myack{}
\fi

{

}

\ifinforms
\section*{Acknowledgments}
\mythanks
\fi

\ifinforms\else

\ifinforms
\ECSwitch
\else
\appendix
\fi
\section{Proofs of Section~\ref{sec:shared_cluster} and Section~\ref{sec:regret} results}
\label{app:regret}

\ifinforms\begin{proof}{Proof of Lemma~\ref{lem:cluster_converge}}\else\begin{proof}[Proof of Lemma~\ref{lem:cluster_converge}]\fi
W.l.o.g., let $K_\tau = K$.
By Assumption~\ref{ass:support} and the Strong Law of Large Numbers (SLLN)~\citep{slln}, for each cluster $k \leq K$,
    $\prob^\infty\left\{ \lim_{t\to\infty}(1/(n_t)) \sum_{i=1}^{n_t}\mathbf{1}\{\hat{u}^i \in S^k_\tau\} = \theta_\star^k\right\}=1. $
     This shows that the cluster weights converge almost surely to $\prob(S^k_\tau) = \theta_\star^k$, the true weights on each cluster support.

Next, as $\prob(S^k_\tau) >0$ for all clusters $k$, we note that $n^k_t \to \infty$ when $t \to \infty$.
Furthermore, since $\prob$ has finite $\pwass$-th moments, the conditional distributions of $\prob$ on these supports $S^k_\tau$ must also have finite $\pwass$-th moments. We denote these conditional distributions $\prob^k$, and their means $\bar{u}^k_\star$. Then, by the SLLN,
$\prob^\infty\left\{ \lim_{n_t^k\to\infty} (1/(n_t^k)) \sum_{\hat{u} \in C^k_t}\hat{u} = \bar{u}^k_\star\right\} = 1,$
 \ie, the cluster means converge almost surely to the true (finite) means of the conditional distributions $\prob^k$.

 Since both the probabilities and the means converge almost surely to finite limits, it follows that
 $\prob^\infty\left\{ \lim_{t\to\infty} \hat{\prob}^K_t =\lim_{t\to\infty}\sum_{k=1}^K\theta^k_t\delta_{\bar{u}_t^k} = \sum_{k=1}^K\theta^k_\star\delta_{\bar{u}_\star^k} =\prob^K_\star\right\} = 1, $
 where $\prob^K_\star$ is constructed to be a discretized distribution of $\prob$ onto $K$ atoms.
 Since $\prob^K_\star$ has finite $\pwass$-th moments by construction, this implies almost sure convergence with respect to the $W_\pwass$ metric~\citep[Theorem 6.9]{optimal_transport_wass}.
\end{proof}

\ifinforms\begin{proof}{Proof of Lemma~\ref{lem:time_diff}}\else\begin{proof}[Proof of Lemma~\ref{lem:time_diff}]\fi
Fix $x \in \mathcal{X}$, and suppress the dependency of $f$ on $x$, \ie,  $f(u)= f(u,x)$. 
We decompose the difference using the intermediate ball $\mathcal{B}'$, which shares its radius $\epsnom_2$ with $\mathcal{B}_2$ and its center $\mathbf{P}_1$ with $\mathcal{B}_1$:
$\max_{\mathbf{Q} \in \mathcal{B}_1}\mathbf{E}_{\mathbf{Q}}[f(u)] - \max_{\mathbf{Q} \in \mathcal{B}_2}\mathbf{E}_{\mathbf{Q}}[f(u)] = \left(\max_{\mathbf{Q} \in \mathcal{B}_1}\mathbf{E}_{\mathbf{Q}}[f(u)] - \max_{\mathbf{Q} \in \mathcal{B}'}\mathbf{E}_{\mathbf{Q}}[f(u)]\right) +\left(\max_{\mathbf{Q} \in \mathcal{B}'}\mathbf{E}_{\mathbf{Q}}[f(u)] - \max_{\mathbf{Q} \in \mathcal{B}_2}\mathbf{E}_{\mathbf{Q}}[f(u)]\right).$
The second term is the fixed-radius discrepancy from~\citep[Theorem 14]{online_dro_s}, applied to $f(u)$. The term is thus bounded by $L\,W_1(\mathbf{P}_1,\mathbf{P}_2) = F_1 W_1(\mathbf{P}_1,\mathbf{P}_2)$ for $p=1$, and $F_p\,W_\pwass(\mathbf{P}_1,\mathbf{P}_2)$ for $p > 1$, where we use the maximum over $\mathcal{X}$ to remove the cited work's dependency on $x$.

We then turn to the first bracketed term. If $\epsnom_1 \leq \epsnom_2$, then $\mathcal{B}_1 \subseteq \mathcal{B}'$, and the difference is nonpositive.  It follows that $L(\epsnom_1-\epsnom_2)_+=0$ is a valid bound. 
Otherwise, suppose $\epsnom_1 > \epsnom_2$. For any arbitrary $\mathbf{Q}_1 \in \mathcal{B}_1$,
we need to construct a distribution $\mathbf{Q}_2 \in \mathcal{B}'$ such that $\mathbf{E}_{\mathbf{Q}_1}[f(u)] - \mathbf{E}_{\mathbf{Q}_2}[f(u)] \leq L(\epsnom_1 - \epsnom_2)$.
To this end, let $\pi \in \Pi(\mathbf{Q}_1,\mathbf{P}_1)$ be an optimal coupling, so that $(v_1,\zeta)\sim\pi$, with $v_1\sim\mathbf{Q}_1$ and $\zeta\sim\mathbf{P}_1$, satisfies
$\big(\mathbf{E}_{\pi}[\|v_1 - \zeta\|^\pwass]\big)^{1/\pwass} = W_\pwass(\mathbf{Q}_1,\mathbf{P}_1) \leq \epsnom_1$.
Next, set $r = \epsnom_2/\epsnom_1 \in (0,1)$, construct $v_2 = \zeta + r(v_1 - \zeta)$, and let $\mathbf{Q}_2$ be the law of $v_2$. Note that by the convexity of $S$, $v_2 \in S$. 
Since $v_2$ is a function of $(v_1,\zeta)$, the pairs $(v_2,\zeta)$ and $(v_1,v_2)$ are, under $\pi$, couplings of $(\mathbf{Q}_2,\mathbf{P}_1)$ and $(\mathbf{Q}_1,\mathbf{Q}_2)$, respectively.
Observe that $W_\pwass(\mathbf{Q}_2,\mathbf{P}_1)^\pwass \leq \mathbf{E}_{\pi}\big[\|v_2 - \zeta\|^\pwass\big] = r^\pwass\,\mathbf{E}_\pi\big[\|v_1 - \zeta\|^\pwass\big] \leq r^\pwass \epsnom_1^\pwass =  (\epsnom_2/\epsnom_1)^\pwass \epsnom_1^\pwass =\epsnom_2^\pwass,$
so $\mathbf{Q}_2 \in \mathcal{B}'$.
Furthermore, by the Lipschitzness of $f$ in $u$ and Jensen's inequality,
\begin{align*}
\mathbf{E}_{\mathbf{Q}_1}[f(u)] - \mathbf{E}_{\mathbf{Q}_2}[f(u)] &= \mathbf{E}_\pi\big[f(v_1)-f(v_2)\big] \leq L\, \mathbf{E}_\pi\big[\|v_1-v_2\|\big] = L(1-r)\,\mathbf{E}_\pi\big[\|v_1-\zeta\|\big]\\
&\leq L(1-r)\big(\mathbf{E}_\pi[\|v_1-\zeta\|^\pwass]\big)^{1/\pwass} \leq L(1-r)\epsnom_1 = L(\epsnom_1-\epsnom_2),
\end{align*}
so $\mathbf{Q}_2$ is constructed as desired. 
Since $\mathbf{Q}_2$ is feasible for $\mathcal{B}'$, $\mathbf{E}_{\mathbf{Q}_2}[f(u)] \leq \max_{\mathbf{Q} \in \mathcal{B}'}\mathbf{E}_{\mathbf{Q}}[f(u)]$, so
$\mathbf{E}_{\mathbf{Q}_1}[f(u)] \leq \max_{\mathbf{Q} \in \mathcal{B}'}\mathbf{E}_{\mathbf{Q}}[f(u)] + L(\epsnom_1-\epsnom_2).$
Taking the supremum over $\mathbf{Q}_1 \in \mathcal{B}_1$ gives the result. 
\end{proof}

\ifinforms\begin{proof}{Proof of Lemma~\ref{thm:nominal_performance}}\else\begin{proof}[Proof of Lemma~\ref{thm:nominal_performance}]\fi
Let $H^{K,\star}_t = \min_{x \in \mathcal{X}} \max_{\mathbf{Q} \in \mathcal{P}^K_{t}} \mathbf{E}_{\mathbf{Q}}[f(\unc, x)]$ denote the exact optimal value of the compressed problem at time $t$. Since $x^K_t$ is only a $\delta$-optimal solution, the achieved value satisfies $H^K_t = \max_{\mathbf{Q} \in \mathcal{P}^K_{t}} \mathbf{E}_{\mathbf{Q}}[f(\unc, x^K_t)] \leq H^{K,\star}_t + \delta$, while the nominal solution $x_t$ is exact. Then,
\begin{equation*}
    H_t = \max_{\mathbf{Q} \in \mathcal{P}_{t}} \mathbf{E}_{\mathbf{Q}}[f( \unc,x_t)] \leq \max_{\mathbf{Q} \in \mathcal{P}_{t}} \mathbf{E}_{\mathbf{Q}}[f(\unc, x^K_{t})] \leq \max_{\mathbf{Q} \in \mathcal{P}^K_{t}} \mathbf{E}_{\mathbf{Q}}[f(\unc, x^K_{t})]  + e_{t} = H^K_t + e_t.
\end{equation*}
The first inequality follows from the optimality of $x_t$, the second follows from~\citep[Theorem 5]{wang2022mean}, 
and $e_t$ is defined in Definition~\ref{def:cluster_value}.
Since $f$ is Lipschitz, using Kantorovich-Rubinstein (KR) duality~\citep{kantorovich_duality}, Lemma~\ref{lem:time_diff}, and Definition~\ref{def:gap}, we also note
$ \max_{\mathbf{Q} \in \mathcal{P}_{t}} \mathbf{E}_{\mathbf{Q}}[f(\unc, x^K_{t})] \leq  \max_{\mathbf{Q} \in \mathcal{P}^K_{t}} \mathbf{E}_{\mathbf{Q}}[f(\unc, x^K_{t})] + \gap(\mathcal{P}_t,\mathcal{P}^K_t).$
Combining all gives $H_t - H^K_t \leq \underline{\psi}_t$. For the other side, using $H^K_t \leq H^{K,\star}_t + \delta$ and \citep[Theorem 5]{wang2022mean},
\begin{equation*}
\begin{aligned}
    H^{K,\star}_t \leq \max_{\mathbf{Q} \in \mathcal{P}^K_{t}} \mathbf{E}_{\mathbf{Q}}[f(\unc, x_{t})] &\leq  H_t + (\tilde{H}_t - H_t) + {\max_{j \leq J}\frac{M_j}{2}}(c_{t,2})^2,
\end{aligned}
\end{equation*}
where $\tilde{H}_t$ is the DRO value with unbounded support. By Lemma~\ref{lem:time_diff}, noting that the bounded support restricts the uncertainty set to an effective radius $\epsnom$ satisfying $\epsnom_t \geq \epsnom \geq 0$, $\tilde{H}_t - H_t \leq L\epsnom_t$.
Next, KR duality and Lemma~\ref{lem:time_diff} again give the remaining bounds, so that they combine to form $H^{K,\star}_t \leq H_t + \bar\psi_t$.
Combining with $H^K_t \leq H^{K,\star}_t + \delta$ gives $H^K_t - H_t \leq \bar\psi_t + \delta$.  
\end{proof}

\ifinforms\begin{proof}{Proof of Theorem~\ref{lem:regret}}\else\begin{proof}[Proof of Theorem~\ref{lem:regret}]\fi
For brevity, let us denote the following shorthands for any $t$:
 $\bar{f}_t(x) = \max_{\mathbf{Q} \in \mathcal{P}_t} \mathbf{E}_{\mathbf{Q}}[f(\unc,x)]$, $\bar{f}^K_t(x) = \max_{\mathbf{Q} \in \mathcal{P}^K_{t}} \mathbf{E}_{\mathbf{Q}}[f(\unc,x)],$
 and the corresponding solutions
 $x_t = \argmin_{x \in \mathcal{X}}\bar{f}_t(x)$, $x^K_{t} = \argmin_{x \in \mathcal{X}}\bar{f}^K_{t}(x),$
 where $x^K_t$ is a $\delta$-optimal minimizer.
We wish to bound
$R(T,K) = \frac{1}{T}\sum_{t=0}^{T-1} \left(\bar{f}_{t+1}(x^K_{t}) - \bar{f}_{t+1}(x_{t+1})\right),$
which can be rewritten as
\begin{equation*}
\begin{aligned}
  R(T,K) &=\frac{1}{T}\sum_{t=0}^{T-1} \left(\bar{f}_{t+1}(x^K_{t}) - \bar{f}_{t}(x_{t})\right) + \frac{1}{T}\sum_{t=0}^{T-1} \left( \bar{f}_{t}(x_{t}) - \bar{f}_{t+1}(x_{t+1})\right).
\end{aligned}
\end{equation*}
Note that the second summation telescopes to $(1/T)\left( \bar{f}_{0}(x_{0}) - \bar{f}_{T}(x_{T})\right).$
We then begin by bounding the terms in the first summation. Using Lemma~\ref{thm:nominal_performance}, we first observe that
$\bar{f}_{t+1}(x^K_{t}) \leq \bar{f}^K_{t+1}(x^K_{t}) + \underline{\psi}_{t+1}$ and $\bar{f}_{t}(x_{t}) \geq \bar{f}_{t}^K(x_{t}) - \bar{\psi}_t.$
Next, since $x^K_{t}$ is a $\delta$-optimal solution, we have $\bar{f}_{t}^K(x^K_{t}) \leq \min_{x\in\mathcal{X}}\bar{f}_{t}^K(x) + \delta \leq \bar{f}_{t}^K(x_{t}) + \delta$, \ie, $\bar{f}_{t}^K(x_{t}) \geq \bar{f}_{t}^K(x^K_{t}) - \delta$.
Furthermore, by KR duality, Lemma~\ref{lem:time_diff}, and Definition~\ref{def:gap}, $ \bar{f}^K_{t+1}(x^K_{t}) -\bar{f}_{t}^K(x^K_{t})\leq \gap(\mathcal{P}^K_{t+1},\mathcal{P}^K_t)$. 
Combining, we obtain
    \begin{equation}
    \begin{aligned}
    \label{lem:regret_eq}
    \frac{1}{T}& \sum_{t=0}^{T-1} (\bar{f}_{t+1}(x^K_{t}) - \bar{f}_{t}(x_{t})) \leq \frac{1}{T}\sum_{t=0}^{T-1}  (\underline{\psi}_{t+1} + \bar{\psi}_t  + \gap(\mathcal{P}^K_{t+1},\mathcal{P}^K_t)) + \delta.
    \end{aligned}
    \end{equation}
For the final term, we note by KR duality and Lemma~\ref{lem:time_diff},
$\bar{f}_{0}(x_{0}) - \bar{f}_{T}(x_{T}) \leq \bar{f}_{0}(x_{T}) - \bar{f}_{T}(x_{T})\leq \gap(\mathcal{P}_0,\mathcal{P}_T)$.
Combining~\eqref{lem:regret_eq} with this bound gives the result.
\end{proof}

\ifinforms\begin{proof}{Proof of Lemma~\ref{thm:clustering_conv}}\else\begin{proof}[Proof of Lemma~\ref{thm:clustering_conv}]\fi
By the continuity of the Wasserstein-$\pwass$ metric, and by the almost-sure convergences of Theorem~\ref{thm:calc_eps} and Lemma~\ref{lem:cluster_converge}, we obtain $\mathbf{P}^\infty\{ \lim_{t \rightarrow \infty } W_{\pwass}(\hat{\prob}_t,\hat{\prob}^K_t)   = W_{\pwass}(\prob,\prob^K_\star)\} = 1$~\citep[Theorem 7.12]{Villani2003Topics}.
For the lower bound, note that $\prob^K_\star$ is supported on at most $K$ atoms, so the definition of the minimum quantization error gives $w_K(\prob) \leq W_{\pwass}(\prob,\prob^K_\star)$.
Lastly, the asymptotic result for $K \to \infty $ holds from~\citep[Lemma 6.1]{quantization}.
\end{proof}

\ifinforms\begin{proof}{Proof of Lemma~\ref{thm:asy_cluster}}\else\begin{proof}[Proof of Lemma~\ref{thm:asy_cluster}]\fi
By the consistency of our clustering procedure (Lemma~\ref{lem:cluster_converge}), the 0-measure of the cluster boundaries (Assumption~\ref{ass:support}), the continuity of linear functions and $\|\cdot\|^{\pwass}$, and the existence of the accumulation points, the limit (and limit superior) of the functions within the summations of Definition~\ref{def:cluster_value} converge almost surely to the functions within the expectations, for $\prob$-almost every $u$. That is, defining
  $\chi_t(u) =  \sum_{k=1}^{K_t}\ones\{u \in C^k_t\}\| u - \bar{u}^k_t\|_2^{\pwass}$ for the first result, we have that $\lim_{t\to\infty} \chi_t(u) = \|u - \bar{u}^{\kappa(u)}_\star\|^{\pwass}_2 = \chi_\star(u)$, for $\prob$-almost every $u$.
  Next, to show $\lim_{t\to\infty} \mathbf{E}_{\hat{\prob}_t}[\chi_t(u)] = \mathbf{E}_{\prob}[\chi_\star(u)]$, we require the function $\chi_t(u)$ to have an integrable dominating function.
  W.l.o.g., let $u$ belong to the cluster $k$. By the triangle inequality, $|\chi_t(u)| = \|u - \bar{u}^k_t\|^{\pwass}_2 \leq (\|u\|_2 + \|\bar{u}^k_t\|_2 )^{\pwass} \leq \bar{C}(1+\|u\|_2^{\pwass} )$ for some constant $\bar{C}$. By the assumption of finite $\pwass$-th moments, this dominating function is integrable. The result then follows from \citep[Theorem 7.12(iv)]{Villani2003Topics}, which establishes a generalized convergence of expectation. Since we have shown almost sure pointwise convergence in the first step, the cited theorem holds even without the assumption on continuity. 
\end{proof}

\ifinforms\begin{proof}{Proof of Theorem~\ref{cor:regret}}\else\begin{proof}[Proof of Theorem~\ref{cor:regret}]\fi
    The convergence of the clustering discrepancies follows from the previous lemmas. The existence of $\hat{x}^K_\star$ follows from the compactness of $\mathcal{X}$.
Then, by the almost sure convergence of
     $\underline{\psi}_t + \bar{\psi}_t$ to $\underline{\psi}_\star + \bar{\psi}_\star$,
    $W_1(\hat{\prob}_0,\hat{\prob}_T)/T$ to $0$, and $W_1(\hat{\prob}^K_t,\hat{\prob}^K_{t+1})$ to $0$, by the Ces\`aro Mean Theorem~\citep{cesaro},
    $R(T,K)$ converges to some value that is $O(\max\{\underline{\psi}_
    \star,\bar{\psi}_\star\})$ with probability 1.
The result for $K \rightarrow \infty$ follows from Lemma~\ref{thm:clustering_conv}.
\end{proof}

The proof of Corollary~\ref{cor:regret_opt} requires the following lemma.
\begin{lemma}
\label{thm:converge_K}
Let $\eps_{t} \sim O(\epsnom_{t})$, where $\epsnom_{t}$ is computed according to Theorem~\ref{thm:calc_eps}. $\prob^\infty$-almost surely we have $H^K_t \rightarrow H^K_\star$ as $t\to\infty$, where $H^K_\star$ is the optimal value of the stochastic optimization problem~\eqref{eq:opt} with distribution $\prob^K_\star$. Furthermore, if $f(u,x)$ is lower-semicontinuous in $x$ for all $u \in \supp$, then any accumulation point $\hat{x}^K_\star = \lim_{t\to \infty}x_t^K$ is almost surely, with respect to $\prob^\infty$, an optimal solution to the stochastic optimization problem with distribution $\prob^K_\star$.
\end{lemma}

\ifinforms\begin{proof}{Proof of Lemma~\ref{thm:converge_K}}\else\begin{proof}[Proof of Lemma~\ref{thm:converge_K}]\fi
This follows from Lemma~\ref{lem:conditional_Conv}, proved below, and an adaptation of~\citep[Theorem 3.6]{mohajerin2018data}, by treating $\prob^K_\star$ as the true distribution.
\end{proof}

\ifinforms\begin{proof}{Proof of Corollary~\ref{cor:regret_opt}}\else\begin{proof}[Proof of Corollary~\ref{cor:regret_opt}]\fi
For the first result, when $f$ is maximum-of-affine and $\supp = \reals^d$, known conjugation results show that $z^{jk} = P^{jk}x$ for constant matrices $P^{jk}$~\citep[Corollary 5.1]{mohajerin2018data}.
Then the accumulation points $\hat{z}^{jk}_\star$ exist as continuous functions of an accumulation point $\hat{x}^K_\star$, which exist by Lemma~\ref{thm:converge_K}.
For the second result, note that $e_t$ is computed with respect to the clustering-induced coupling, \ie, each $\hat{u}$ is associated with the mean $\bar{u}^k_t$ of the cluster it belongs to. If this is equivalent to the optimal coupling, then by H\"{o}lder's inequality, for each $t, j$, and $\hat{u}$,
$(-z^{jk}_t)^T (\hat{u} - \bar{u}^k_t) \leq \|z^{jk}_t\|_*\|\hat{u}- \bar{u}^k_t\|\leq L_j\|\hat{u} - \bar{u}^k_t\| =  L_j \min \left \{\|\hat{u} - \bar{u}\|\mid \bar{u} \in \{\bar{u}^k_t\}_{k=1}^K \right\}$
where $\|z^{jk}_t\|_*\leq L_j$~$\forall j\leq J$, $k\leq K$, and $t\geq1$ by the conjugates of max-of-affine functions and the compactness of $\mathcal{X}$. This implies
$e_t \leq L W_{1}(\hat{\prob}_t,\hat{\prob}^K_t)$. 
The result follows by taking the limsup.
\end{proof}

\ifinforms\begin{proof}{Proof of Lemma~\ref{lem:conditional_Conv}}\else\begin{proof}[Proof of Lemma~\ref{lem:conditional_Conv}]\fi
    By assumption, all cluster supports have nonzero measure, \ie, $\prob(S^k_\tau) > 0$.
    Since $\prob$ is light-tailed, the constant $R$ of Assumption~\ref{ass:light_tail} is finite for some $\alpha > 1$.
    It follows that the conditional distributions of $\prob$ on the supports $S^k_\tau$, denoted $\prob^k$, are also light-tailed, with the same $\alpha$ and a constant $R^k \leq O(R)$, \ie,
    $R^k = \int_{S_\tau^k} \exp(\|u\|^{\alpha})\prob(u|u\in S^k_\tau)(du) < \infty.$
    By Theorem~\ref{thm:calc_eps},
    $\prob^{t}(W_{\pwass}(\prob^k,\hat{\prob}^{t,k})\leq O(\epsnom_t)) \geq 1-\beta_t,$
    where $\hat{\prob}^{t,k}$ is the conditional distribution of $\hat{\prob}_t$ on support $S^k_\tau$.
    By KR duality and the ordering of Wasserstein distances,
    $\left|\Expect_{\prob^k}[u_i] - \Expect_{\hat{\prob}^{t,k}}[u_i] \right| \leq W_1(\prob^k,\hat{\prob}^{t,k}) \leq W_{\pwass}(\prob^k,\hat{\prob}^{t,k}),$ for $i =1,\dots,d.$
    This implies that the means $\bar{u}^k_t$ converge to the true means $\bar{u}^k_\star$ with rate $O(\epsnom_t)$, with respect to the infinity norm, \ie,
    $\| \bar{u}^k_t - \bar{u}^k_\star\|_\infty \leq O(\epsnom_t).$
    Then, by the equivalence of norms on $\reals^d$ and the triangle inequality through the intermediate distribution $\sum_{k=1}^K\theta^k_t\delta_{\bar{u}^k_\star}$,
    $$W_{\pwass}(\hat{\prob}^K_t,\prob^K_\star) \leq \left(\sum_{k=1}^K\theta^k_t \| \bar{u}^k_t - \bar{u}^k_\star\|^{\pwass}\right)^{1/{\pwass}} + W_{\pwass}\Big(\sum_{k=1}^K\theta^k_t\delta_{\bar{u}^k_\star},\,\sum_{k=1}^K\theta^k_\star\delta_{\bar{u}^k_\star}\Big).$$
    The first term transports each atom $\bar{u}^k_t$ to $\bar{u}^k_\star$ at unchanged weights, and is $O(\epsnom_t)$.
    The second term only redistributes weight among the $K$ fixed atoms $\bar{u}^k_\star$, so it is at most $\bar{B}\big(\sum_{k=1}^K|\theta^k_t-\theta^k_\star|\big)^{1/\pwass}$, where $\bar{B}$ is the diameter of the atom set.
    Since each $\theta^k_t$ is an empirical frequency, Hoeffding's inequality and a union bound over the $K$ clusters give $\sum_{k=1}^K|\theta^k_t-\theta^k_\star| \leq O(\sqrt{\log(\beta_t^{-1})/t}\,)$ with probability at least $1-\beta_t$, so the second term is $O((\log(\beta_t^{-1})/t)^{1/(2\pwass)})$.
\end{proof}

\ifinforms\begin{proof}{Proof of Theorem~\ref{thm:regret}}\else\begin{proof}[Proof of Theorem~\ref{thm:regret}]\fi
We fix $T \gg \tau$, and bound the terms of $R(T,K)$ using the $W_1$ terms in $\underline{\psi}_{t}$ and $\bar{\psi}_t$ as upper bounds.
By assumption, $\epsnom_t\sim (\log(\beta_t^{-1})/t)^{\min(\pwass/d,1/2)}$, so their averages are $O((\log(\beta_T^{-1})/T)^{\min(\pwass/d,1/2)})$, and the summability of $\beta_t$ makes its decay at least sublinear.
Every remaining term is a Wasserstein distance between two empirical objects, bounded by inserting $\prob$ or $\prob^K_\star$ and applying the triangle inequality. Theorem~\ref{thm:calc_eps} gives $W_1(\hat{\prob}_t,\prob)\leq\epsnom_t$ with probability $1-\beta_t$, and Lemma~\ref{lem:conditional_Conv} gives $W_1(\prob^K_\star,\hat{\prob}^K_t)\leq\tilde{\epsnom}_t\leq O(\epsnom_t)$ with probability $1-\beta_t$ for $t\geq\tau$. 
Therefore, the terms in $W_{1}(\hat{\prob}_t,\hat{\prob}^K_t)$ and $W_{1}(\hat{\prob}_{t-1},\hat{\prob}^K_{t-1})$ are together $O(W_1(\prob,\prob^K_\star) + (\log(\beta_{T}^{-1})/T)^{\min(\pwass/d,1/2)})$. Next, splitting $W_1(\hat{\prob}_0,\hat{\prob}_T)$ through $\prob$ leaves one term constant in $T$ and one bounded by $\epsnom_T$, so its average is $O(T^{-1})$, and splitting $W_1(\hat{\prob}^K_t,\hat{\prob}^K_{t+1})$ through $\prob^K_\star$ bounds it by $2\tilde{\epsnom}_{t}$, averaging to $O(\epsnom_T)$.
Intersecting all events, $R(T,K) \leq O((\log(\beta_T^{-1})/T)^{\min(\pwass/d,1/2)})+ O(W_1(\prob,\prob^K_\star))$ with probability $1-\bar{\beta} \geq \Pi_{t=1}^T(1-\beta_t)^4$, or $1-\sum_{t=1}^T 4\beta_t$ by the Bonferroni Inequality~\citep{Bonferroni}.
\end{proof}

\section{Proofs of Section~\ref{sec:subgrad_sec} results}
\label{app:subgrad}

\ifinforms\begin{proof}{Proof of Proposition~\ref{thm:gen_subgrad}}\else\begin{proof}[Proof of Proposition~\ref{thm:gen_subgrad}]\fi
Let $j$ index an active piece of cluster $k$'s oracle at its maximizer $\zeta^k$.
Any other piece $j'$ satisfies $a_{j'}(\bar x)^T\zeta^{k}+b_{j'}(\bar x)\leq a_{j}(\bar x)^T\zeta^{k}+b_{j}(\bar x)$, and since both incur the same additional penalty at $\zeta^{k}$, piece $j$ also attains the maximum $f(\zeta^{k},\bar x)$.
The function $x \mapsto a_{j}(x)^T\zeta^{k}+b_{j}(x)$ is therefore a convex minorant of $f(\zeta^{k},\cdot)$ that coincides with it at $\bar x$, and its subgradients at $\bar x$ belong to $\partial_x f(\zeta^{k},\bar x)$.
For $p=\infty$, Danskin's theorem~\citep{danskin1967theory} then gives~\eqref{eq:general_subgrad} directly.
For finite $p$, the chain rule adds to~\eqref{eq:general_subgrad} the cross term $\partial_\lambda h(\bar{x},\lambda^\star(\bar{x}))\cdot\partial_x\lambda^\star(\bar{x})$, which we show is $0$ whenever $p>1$ or $S\subsetneq\reals^d$.
In those cases, $\varphi_\lambda$ is finite for $\lambda>0$ (see the proof of Proposition~\ref{thm:finitep}), and $h(\bar{x},\lambda)$ is finite and convex, so either $\lambda^\star(\bar x)>0$ is an interior stationary point, or $\lambda^\star(\bar x)=0$.
If $\lambda^\star(\bar{x})>0$ is an interior stationary point, \ie, $0\in\partial_\lambda h(\bar{x},\lambda^\star(\bar{x}))$, then by the envelope theorem for the infimum over $\lambda$, the cross term is $0$; if $\lambda^\star(\bar{x})=0$, the cross term is also $0$.
\end{proof}

\ifinforms\begin{proof}{Proof of Lemma~\ref{thm:combined_tol}}\else\begin{proof}[Proof of Lemma~\ref{thm:combined_tol}]\fi
Let $H(x) \define \max_{\mathbf{Q}\in\mathcal{P}^K}\mathbf{E}_{\mathbf{Q}}[f(u,x)] = \inf_{\lambda\geq0}h(x,\lambda)$, and $H^\star \define H(\bar x) = h(\bar{x},\lambda^\star(\bar x))$.
Since $\partial_\lambda h(x,\lambda) \in [\epsnom^p-R_S^p,\epsnom^p]$, which is bounded by $ R_S^p \geq \epsnom^p$,  $h$ is Lipschitz in $\lambda$.

We begin with the value tolerance. Since $|\hat\lambda-\lambda^\star(\bar x)|\leq\sigma$, Lipschitzness gives $H^\star \leq h(\bar x,\hat\lambda) \leq H^\star+R_S^p\sigma$.
On the other hand, the value-tolerance $\Gamma$ acts on the inner maximization, so the approximate solution is {\it smaller} than the optimal. The bound takes the maximum, $\max\{\Gamma,R_S^p\sigma\}$.

For the subgradient result, let $\hat\varphi(x) \define \sum_{k=1}^K\theta^k f(\hat\zeta^k,x)$, so that $\hat g\in\partial\hat\varphi(\bar x)$ and $\hat\varphi(x)\geq\hat\varphi(\bar x)+\hat g^T(x-\bar x)$. Let $c \define \hat\lambda\sum_k\theta^k\|\hat\zeta^k-\bar u^k\|^p$.  By the value-tolerance $\Gamma$, $\hat\lambda\epsnom^p+\hat\varphi(\bar x) \geq H^\star-\Gamma +c$, and by the feasibility of $\hat\zeta^k$,  $\hat\lambda\epsnom^p + \hat\varphi(x) \leq h(x,\hat\lambda)+c$. Subtracting the expressions, we obtain $\hat\varphi(x)-\hat\varphi(\bar x) \leq h(x,\hat\lambda)-H^\star+\Gamma$.
By the Lipschitzness of $h$ in $\lambda$, $h(x,\hat\lambda) \leq h(x,\lambda^\star(x))+R_S^p|\hat\lambda-\lambda^\star(x)| = H(x)+R_S^p|\hat\lambda-\lambda^\star(x)|$.
By the triangle inequality, $|\hat\lambda-\lambda^\star(x)| \leq |\hat\lambda-\lambda^\star(\bar x)|+|\lambda^\star(\bar x)-\lambda^\star(x)| \leq \sigma+|\lambda^\star(\bar x)-\lambda^\star(x)|$.
Combining, $h(x,\hat\lambda) \leq H(x)+R_S^p(\sigma+|\lambda^\star(\bar x)-\lambda^\star(x)|)$, so
$\hat\varphi(x)-\hat\varphi(\bar x) \;\leq\; H(x)-H^\star+\delta_{\text{grad}}.$
Combining with $\hat g^T(x-\bar x) \leq \hat\varphi(x)-\hat\varphi(\bar x)$ gives~\eqref{eq:approx_subgrad_ineq}. Lastly, when $p=\infty$, $R_S^p\sigma=0$, and the result follows.
\end{proof}

\ifinforms\begin{proof}{Proof of Proposition~\ref{thm:finitep}}\else\begin{proof}[Proof of Proposition~\ref{thm:finitep}]\fi
For $S=\reals^d$, substituting $\zeta=u+v$, setting $r = \|v\|$, and using $\sup_{\|v\|=r}a_j(\bar x)^Tv=r\|a_j(\bar x)\|_*$, gives $\varphi_\lambda^j(u)=a_j(\bar x)^Tu+b_j(\bar x)+\sup_{r\geq0}(r\|a_j(\bar x)\|_*-\lambda r^p)$.
At $p=1$, $r\mapsto r(\|a_j(\bar x)\|_*-\lambda)$ is linear, so this supremum is $0$ if $\|a_j(\bar x)\|_*\leq\lambda$ and $+\infty$ otherwise. At $p>1$, it is concave in $r$, maximized in closed form at $\frac1q\|a_j(\bar x)\|_*^q/(\lambda p)^{q-1}$. The cost is $O(1)$.
For $S\subsetneq\reals^d$, the dual depends on the specific support; see Appendix~\ref{sec:inner_derivations}.
In all cases except $p=1$ with $S=\reals^d$, $\varphi_\lambda$ is finite for $\lambda>0$, and $h(\bar{x},\lambda)$ is finite and convex, with monotone nondecreasing subgradient ${\epsnom}^p-\sum_k\theta^k\|\zeta^{k}(\lambda)-\bar u^k\|^p$.
It follows that either $\hat{\lambda}$ is an interior stationary point of $h(\bar{x},\lambda)$, located by bisecting on the sign of its subgradient, or $\lambda^\star=0$ (if $S\subsetneq\reals^d$ and the subgradient of $h$ at $\lambda=0$ is already nonnegative).
Therefore, with optimality tolerance $\sigma$ in the bisection search, $T_{\rm outer}=O(\log(1/\sigma))$.
For $p=\infty$, the closed form at $S=\reals^d$ follows from $\sup_{\|\zeta-u\|\leq\epsnom}a_j(\bar x)^T\zeta = a_j(\bar x)^Tu+\epsnom\|a_j(\bar x)\|_*$, and since the Wasserstein-$\infty$ ball decouples across pieces and clusters, no outer search for a shared $\lambda$ is required.
Note that for fixed $\lambda$, the inner solves are independent of one another, and could be parallelized.
\end{proof}

\ifinforms\begin{proof}{Proofs of Corollaries~\ref{cor:p1_closed} and~\ref{cor:pinf_closed}}\else\begin{proof}[Proofs of Corollaries~\ref{cor:p1_closed} and~\ref{cor:pinf_closed}]\fi
By the proof of Proposition~\ref{thm:finitep}, $h(\lambda)$ is finite only for $\lambda\geq \max_{j\leq J}\|a_j(\bar x)\|_*$,  where it equals the strictly increasing function $\lambda\epsnom+\sum_k\theta^kf(\bar u^k,\bar x)$. Therefore, $\lambda^\star=\max_{j\leq J}\|a_j(\bar x)\|_*$, and the first result follows. The $p=\infty$ result follows from $\varphi_\infty^j(\bar u^k)=a_j(\bar x)^T\bar u^k+b_j(\bar x)+\epsnom\|a_j(\bar x)\|_*$ at $S=\reals^d$.
\end{proof}

\ifinforms\begin{proof}{Proof of Theorem~\ref{thm:regret_single_resolve}}\else\begin{proof}[Proof of Theorem~\ref{thm:regret_single_resolve}]\fi
Let $\bar{f}_t(x) = \max_{\mathbf{Q} \in \mathcal{P}_t}\mathbf{E}_{\mathbf{Q}}[f(\unc,x)]$ and $\bar{f}^K_t(x) = \max_{\mathbf{Q} \in \mathcal{P}^K_{t}}\mathbf{E}_{\mathbf{Q}}[f(\unc,x)]$, and denote
$H^{K,\star}_t = \min_{x \in \mathcal{X}} \bar{f}^K_t(x),~x^{K,\star}_t \in \argmin_{x \in \mathcal{X}} \bar{f}^K_t(x)$ as the exact optimal value and solution of~\eqref{eq:robustopt_p_max} at time $t$.
Let $x^\star = x^{K,\star}_0$.
By Lemma~\ref{thm:combined_tol}, with $\hat g=g_t$, $\bar x=x^K_{t-1}$, and $x=x^\star$,
$\bar{f}^K_t(x^\star) \;\geq\; \bar{f}^K_t(x^K_{t-1}) + g_t^T(x^\star - x^K_{t-1}) - \delta_{\text{grad},t}.$
\ie, $g_t$ is a $\delta_{\text{grad},t}$-subgradient of $\bar{f}^K_t$ at $x^K_{t-1}$.
By Assumption~\ref{ass:subgrad_lip}, $\|g_t\|_2 \leq L_x$, since $\bar{f}^K_t$ is a pointwise maximum of functions $L_x$-Lipschitz in $x$.
The iterates $x^K_t = \Pi_{\mathcal{X}}(x^K_{t-1}-\eta_tg_t)$ are therefore those of online PGD on the sequence $\bar{f}^K_1,\dots,\bar{f}^K_T$, with step size $\eta_t = R_{\mathcal{X}}/(L_x\sqrt t)$, so the standard regret bound against the fixed comparator $x^\star$~\citep[Theorem 3.1]{hazan2016oco} applies, with diameter $R_{\mathcal{X}}$ and gradient bound $L_x$. We add the $\delta_{\text{grad},t}$ tolerances.
\begin{equation}
\label{eq:sgd_sum}
\sum_{t=1}^{T}\left(\bar{f}^K_t(x^K_{t-1}) - \bar{f}^K_t(x^\star)\right) \leq \frac{3R_{\mathcal{X}}L_x\sqrt T}{2} + \sum_{t=1}^T\delta_{\text{grad},t}.
\end{equation}
Next, by Lemma~\ref{lem:time_diff} and KR duality applied at $x=x^\star$ for every $t \leq T$,
$\bar{f}^K_t(x^\star) - \bar{f}^K_0(x^\star) \leq  \gap(\mathcal{P}^K_t,\mathcal{P}^K_0).$
Since $\bar{f}^K_0(x^\star)= H^{K,\star}_0$, this gives $\bar{f}^K_t(x^\star) \leq H^{K,\star}_0 +  \gap(\mathcal{P}^K_t,\mathcal{P}^K_0)$.
By the same argument applied at $x=x^{K,\star}_t$, and using $\bar{f}^K_0(x^{K,\star}_t) \geq H^{K,\star}_0$,
$$H^{K,\star}_t = \bar{f}^K_t(x^{K,\star}_t) \geq \bar{f}^K_0(x^{K,\star}_t) -  \gap(\mathcal{P}^K_0,\mathcal{P}^K_t)\geq H^{K,\star}_0 -  \gap(\mathcal{P}^K_0,\mathcal{P}^K_t),$$
where we note that the $\gap(\cdot)$ indices switch order.
Combining these expressions, $\bar{f}^K_t(x^\star) - H^{K,\star}_t \leq \gap(\mathcal{P}^K_t,\mathcal{P}^K_0) + \gap(\mathcal{P}^K_0,\mathcal{P}^K_t) \leq 2 \gap(\mathcal{P}^K_0,\mathcal{P}^K_t)$ for every $t \leq T$, using $\gap(\mathcal{P}^K_t,\mathcal{P}^K_0) \leq \gap(\mathcal{P}^K_0,\mathcal{P}^K_t)$ as $\epsnom_t \leq \epsnom_0$. Adding this to~\eqref{eq:sgd_sum} and noting that the LHS cancels out as $\bar{f}^K_t(x^K_{t-1}) - \bar{f}^K_t(x^\star) + \bar{f}^K_t(x^\star) - H^{K,\star}_t = \bar{f}^K_t(x^K_{t-1}) - H^{K,\star}_t$,
\begin{equation}
\label{eq:sgd_sum_H}
\sum_{t=1}^{T}\left(\bar{f}^K_t(x^K_{t-1}) - H^{K,\star}_t\right) \leq \frac{3R_{\mathcal{X}}L_x\sqrt T}{2} + \sum_{t=1}^T\delta_{\text{grad},t} + 2\sum_{t=1}^T \gap(\mathcal{P}^K_0,\mathcal{P}^K_t).
\end{equation}
We now connect this bound to $R_{\text{approx}}(T,K)$, which evaluates $x^K_t$ against the non-compressed ambiguity sets.
By Lemma~\ref{thm:nominal_performance}, for every $t$, $\bar{f}_t(x^K_{t-1}) \leq \bar{f}^K_t(x^K_{t-1}) + \underline\psi_t$, and $\bar{f}_{t-1}(x_{t-1}) \geq \bar{f}^K_{t-1}(x_{t-1}) - \bar\psi_{t-1} \geq H^{K,\star}_{t-1} - \bar\psi_{t-1}$, using that $H^{K,\star}_{t-1}$ is a global minimum of $\bar{f}^K_{t-1}$.
Subtracting,
$\bar{f}_t(x^K_{t-1}) - \bar{f}_{t-1}(x_{t-1}) \;\leq\; \underline\psi_t + \bar\psi_{t-1} + \left(\bar{f}^K_t(x^K_{t-1}) - H^{K,\star}_{t-1}\right).$
Summing over $T$ and using the telescoping identity $\sum_{t=1}^T(H^{K,\star}_t-H^{K,\star}_{t-1}) = H^{K,\star}_T-H^{K,\star}_0$, we decompose 
$\sum_{t=1}^{T}\left(\bar{f}^K_t(x^K_{t-1}) - H^{K,\star}_{t-1}\right) = \sum_{t=1}^{T}\left(\bar{f}^K_t(x^K_{t-1}) - H^{K,\star}_t\right) + \left(H^{K,\star}_T - H^{K,\star}_0\right).$
The first term on the right is bounded above by~\eqref{eq:sgd_sum_H}.
For the second, since $H^{K,\star}_T = \min_x\bar{f}^K_T(x) \leq \bar{f}^K_T(x^\star)$, and $\bar{f}^K_T(x^\star) \leq H^{K,\star}_0+ \gap(\mathcal{P}^K_T,\mathcal{P}^K_0)$ was already given above (set $t=T$), $H^{K,\star}_T-H^{K,\star}_0 \leq  \gap(\mathcal{P}^K_T,\mathcal{P}^K_0)$.
Combining,
$\sum_{t=1}^{T}\left(\bar{f}_t(x^K_{t-1}) - \bar{f}_{t-1}(x_{t-1})\right) \;\leq\; \sum_{t=1}^{T}\left(\underline\psi_t+\bar\psi_{t-1}\right) + \frac{3R_{\mathcal{X}}L_x\sqrt T}{2} + \sum_{t=1}^T\delta_{\text{grad},t} + 2\sum_{t=1}^T \gap(\mathcal{P}^K_0,\mathcal{P}^K_t)+  \gap(\mathcal{P}^K_T,\mathcal{P}^K_0).$
Adding the telescoped term $\bar{f}_0(x_0)-\bar{f}_T(x_T)$, bounded as in the proof of Theorem~\ref{lem:regret}, and dividing by $T$ yields an upper bound. Note that some terms appear within~\eqref{eq:regret_det}, so substituting in $\tilde{R}(T,K)$ yields~\eqref{eq:regret_single_resolve}. The additional term that appears in $\tilde{R}$ but not in $R_{\text{approx}}$ is nonnegative, so the upper bound holds. 
\end{proof}

\section{Inner-problem derivations}
\label{sec:inner_derivations}
In this section, we derive the entries of Table~\ref{tab:inner_cost}: the per-piece oracle $\varphi_\lambda^j(u)=\sup_{\zeta\in S}(a_j(\bar x)^T\zeta+b_j(\bar x)-\lambda\|\zeta-u\|^p)$ for finite $p$, and $\varphi_\infty^j(u)=\sup_{\zeta\in S,\,\|\zeta-u\|\leq\epsnom}(a_j(\bar x)^T\zeta+b_j(\bar x))$ for $p=\infty$, sectioned by support $S$ and ground norm $\|\cdot\|$.
For simplicity, we write $\beta=a_j(\bar x)$ and drop the additive constant $b_j(\bar x)$, which carries through unchanged.

\paragraph{Unbounded support.}
See the proof of Proposition~\ref{thm:finitep}. 

\paragraph{Box, $\|\cdot\|_1$.}
Let $S=[l,h]$, and $v=\zeta-u$.
Let $\delta_i^+=h_i-u_i$ and $\delta_i^-=u_i-l_i$ denote the box slacks.

\textit{Case $p=1$.} The penalty $\lambda\|v\|_1=\lambda\sum_i|v_i|$ separates additively by coordinate, so we can solve individually $\max_{v_i\in[-\delta_i^-,\delta_i^+]}(\beta_iv_i-\lambda|v_i|)$, giving $v_i^\star=\delta_i^+$ if $\beta_i>\lambda$, $v_i^\star=-\delta_i^-$ if $\beta_i<-\lambda$, and $v_i^\star=0$ otherwise. Cost: $O(d)$.

\textit{Case $p>1$.} Denote $r=\|v\|_1\geq0$, and split the problem in two stages: fix $r$ and solve the linear problem $F(r)\define\max\{\beta^Tv : \|v\|_1\leq r,\ -\delta_i^- \leq v_i \leq \delta_i^+~ \forall i\}$, then maximize $F(r)-\lambda r^p$ over $r$.
The inner problem is the continuous knapsack problem~\citep{dantzig1957discrete}, which can be solved via sorting. 
We sort the coordinates by $|\beta_i|$ in decreasing order, $\pi(1),\dots,\pi(d)$, and let $\rho_0=0$, $\rho_k=\sum_{l\leq k}\delta_{\pi(l)}$ ($\delta_i\define\delta_i^+$ if $\beta_i>0$, $\delta_i^-$ if $\beta_i<0$) be the cumulative slack of the $k$ largest coordinates. For $r\in[\rho_{k-1},\rho_k]$, the greedy allocation gives
$F(r) = \sum_{l<k}|\beta_{\pi(l)}|\delta_{\pi(l)} + |\beta_{\pi(k)}|\,(r-\rho_{k-1})$, a piecewise-linear function of $r$. Maximizing $F(r)-\lambda r^p$ is then a segment-by-segment search in $O(d)$, where on each segment the optimal $r$ has closed form  $(|\beta_{\pi(k)}|/(\lambda p))^{1/(p-1)}$. The total cost is $O(d\log d)$.

\textit{Case $p=\infty$.} Same as above, but with the budget fixed at $r = \epsnom$. Cost: $O(d\log d)$.

\paragraph{Box, $\|\cdot\|_\infty$.}
Let $\delta_i=\delta_i^+$ if $\beta_i>0$ and $\delta_i=\delta_i^-$ if $\beta_i<0$.

\textit{Case $p=1$.} At transport radius $\|v\|_\infty\leq r$, the optimal perturbation is $v_i=\mathrm{sign}(\beta_i)\min(r,\delta_i)$. Summing over coordinates, the reformulated objective $F(r)-\lambda r$, with $F(r) = \sum_i|\beta_i|\min(r,\delta_i)$,
is concave and piecewise linear in $r\geq0$, with slope $\sum_i|\beta_i|\ones\{\delta_i>r\}-\lambda$ that decreases as $r$ grows. Sorting the $\delta_i$ and walking $r$ through the resulting breakpoints in increasing order, the slope decreases by $|\beta_i|$ at each $\delta_i$ passed. Therefore, $r^\star$ is the point at which this running slope first reaches zero, found by a single forward scan of the sorted $\delta_i$ in $O(d)$ time after the $O(d\log d)$ sort.

\textit{Case $p>1$.} The same reduction applies with $F(r)-\lambda r^p$ instead of $F(r)-\lambda r$. The result follows similarly as the $\|\cdot\|_1$ case above, where the cost is $O(d\log d)$, dominated by the sort. 

\textit{Case $p=\infty$.} Since $\|v\|_\infty\leq\epsnom$ bounds every coordinate independently, $v_i\in[-\epsnom,\epsnom]$, so $S\cap B_\infty(u,\epsnom)$ is a box. The problem is coordinate-separable, with $v_i^\star=\min(\delta_i^+,\epsnom)$ if $\beta_i>0$, $v_i^\star=-\min(\delta_i^-,\epsnom)$ if $\beta_i<0$. Cost: $O(d)$.

\paragraph{Box, $\|\cdot\|_2$.}
The penalty $\lambda\|\zeta-u\|_2^p$ is smooth, with gradient
$\nabla_\zeta\big(\lambda\|\zeta-u\|_2^p\big) = \lambda p\,\|\zeta-u\|_2^{p-2}(\zeta-u).$
Defining $r\define\|\zeta^\star-u\|_2$ as the radius attained at the optimum and $\mu\define\lambda p\,r^{p-2}$, KKT conditions set $\zeta_i^\star(\mu)=\mathrm{clip}(u_i+\beta_i/\mu,\,l_i,\,h_i)$.

\textit{Case $p=2$.} Here $\mu=\lambda p\,r^{p-2}=2\lambda$ exactly, independent of $r$. Since the penalty $\lambda\|\zeta-u\|_2^2=\lambda\sum_i(\zeta_i-u_i)^2$ is coordinate-separable, differentiating directly gives $\beta_i-2\lambda(\zeta_i-u_i)=0$. It follows that $\zeta^\star=\zeta^\star(2\lambda)$ is a closed-form clip, with cost $O(d)$.

\textit{Case $p\neq2$.} For $p\neq2$, $\mu=\lambda p\,r(\mu)^{p-2}$ depends on the radius $r(\mu)\define\|\zeta^\star(\mu)-u\|_2$ that $\mu$ itself induces (at $p=1$, this collapses to $\mu\,r(\mu)=\lambda$. See Example~\ref{ex:box_l2_p1}). Since $r(\mu)$ is nonincreasing in $\mu$, the monotone fixed-point equation has a unique root, found by bisecting $\mu\mapsto\mu-\lambda p\,r(\mu)^{p-2}$ to tolerance $\sigma_{\text{in}}$ with cost $O(d\log(1/\sigma_{\text{in}}))$.

\textit{Case $p=\infty$.} The target radius is fixed directly, at $r=\epsnom$. The same bisection on $\mu$ applies, now solving $r(\mu)=\epsnom$, at the same cost $O(d\log(1/\sigma_{\text{in}}))$.

\paragraph{Polyhedral support, all norms and $p$.}
Unlike the box support, the coordinates are no longer decoupled, so a direct solve for each inner problem is already optimal. The cost is $O(d(m+d)^{1.5})$ for LPs / power-cone programs, and $O((m+d)^{1.5})$ for the $\ell_2$-norm SOCPs.

\paragraph{Norm-ball support, $\|\cdot\|_2$.}
Let $S=B(c,R)$.
Defining again $\mu=\lambda p\,r^{p-2}$ and $r=\|\zeta^\star-u\|_2$, the stationarity conditions of the oracle give  $\beta=\mu(\zeta^\star-u)+\eta(\zeta^\star-c)$. 
This is linear in $\zeta^\star$ for fixed $(\mu,\eta)$, solving explicitly to
$\zeta^\star(\mu,\eta) = (\beta+\mu u+\eta c)/({\mu+\eta}).$
To solve for $\mu$ and $\eta$, first check whether the ball constraint is inactive ($\eta=0$): this gives $\mu=(\lambda p)^{1/(p-1)}\|\beta\|_2^{(p-2)/(p-1)}$ in closed form. 
If the resulting $\zeta^\star(\mu,0)=u+\beta/\mu$ satisfies $\|\zeta^\star(\mu,0)-c\|_2\leq R$, it is optimal, at cost $O(d)$.
Otherwise, the ball constraint is active, and $(\mu,\eta)$ must be found jointly, subject to $\|\zeta^\star(\mu,\eta)-u\|_2=r=(\mu/(\lambda p))^{1/(p-2)}$ and $\|\zeta^\star(\mu,\eta)-c\|_2=R$. Bisecting on $\eta$, with each query an $O(d)$ evaluation of $\zeta^\star(\mu,\eta)$, gives $\zeta^\star$ to tolerance $\sigma_{\text{in}}$ in $O(d\log(1/\sigma_{\text{in}}))$, for any finite $p$; at $p=2$, $\mu=2\lambda$ is fixed and only $\eta$ is bisected.

\textit{Case $p=\infty$.} Fixing $r=\epsnom$ directly, the same two-multiplier system as above is solved by bisecting on $\eta$ to find the target radius, in $O(d\log(1/\sigma_{\text{in}}))$.

\section{Online clustering algorithm (Section~\ref{sec:clustream})}
\label{app:onlinecluster}
We give details for the online clustering algorithm, introduced in Section~\ref{sec:clustream}.
We keep track of two sets of clusters: a set of $Q \geq K$ {\it microclusters}, and a set of $K$ {\it macroclusters}. 

\paragraph{Initialization}
We initialize by solving approximately the $k$-centers problem~\eqref{def:k_centers} with respect to $\hat{\prob}_0$, with up to $Q$ clusters. This set of centers is denoted $A^Q_0$. We assign points to their closest center, and define the support of the $q$-th cluster $S^q_0$ as the Voronoi region $V(a_0^q)$ around its center $a_0^q \in A^Q_0$.
Each microcluster stores its center, mean, root-mean-squared error (RMSE) $\eta^q_t$, and weight; for singleton clusters, the RMSE is heuristically initialized as twice the minimum stored RMSE.

To create $K$ {\it macroclusters}, we solve approximately the $k$-centers problem with respect to $A^Q_0$, and obtain a set of centers $A^K_0$. Each microcluster is assigned to the macrocluster with the closest center; a macrocluster then contains the points of its constituent microclusters, with their combined weight, and the weighted average of their means. The support of the $k$-th macrocluster is the Voronoi region around the $k$-th center $a_0^k \in A^K_0$, plus the datapoints assigned to it, and minus the other datapoints, \ie,
$ S^k_t = \big(V(a_t^k) \cup C^k_t\big) \setminus \textstyle\bigcup_{k' \neq k} C^{k'}_t.$
With this, the supports satisfy Assumption~\ref{ass:support}. The adjustment is needed due to the two-layer structure: a datapoint clustered into microcluster $q$ and assigned to macrocluster $k$ may be closer to the center of another macrocluster $k' \neq k$. Note that for continuous supports, such adjustments are negligible, so datapoints can be discarded.

\paragraph{Online updating procedure}
For $t \geq 1$, when we observe a new datapoint $\hat{u}$, we compute its distance to the center of the microcluster whose support contains it, \ie, if $\hat{u} \in S^q_{t-1}$,
$d_t = \|\hat{u} - a^q_{t-1}\|.$
If $d_t \leq \rho\,\eta^q_{t-1}$, where $\eta^q_{t-1}$ is the cluster's RMSE and $\rho \geq 1$ is a fixed multiplier ($\rho = 1.25$ in our experiments), we assign $\hat{u}$ to this cluster.
We then update all microcluster weights, as well as the mean and RMSE of the assigned microcluster. The macroclusters are adjusted accordingly, based on their constituent microclusters. 
Otherwise, we create a new microcluster with this datapoint as its center, and assign it the weight $1/n_t$. The weights of the other microclusters are decreased accordingly, as in the Reclustering weight update of Section~\ref{sec:recluster}. If the number of microclusters then exceeds $Q$, we merge the two with the closest centers: the merged center is the weighted average of the two, with mean, RMSE, and weight combined accordingly, and the supports reassigned as above. 
Every 50 steps, we perform approximate $k$-centers on the $Q$ microclusters, warm-starting at the previous centers, and refresh the $K$ macroclusters. 
After time $\tau < \infty$, we only cluster points based on the latest supports $\mathcal{S}^K_\tau$.

 \fi

\end{document}